\documentclass[11pt]{amsart}

\usepackage{amsmath}
\usepackage{amssymb}

\newtheorem{thm}{Theorem}[section]
\newtheorem{cor}[thm]{Corollary}
\newtheorem{lem}[thm]{Lemma}
\newtheorem{prop}[thm]{Proposition}

\theoremstyle{definition}
\newtheorem{defn}[thm]{Definition}

\newtheorem{ex}[thm]{Example}

\theoremstyle{remark}
\newtheorem{rem}[thm]{Remark}

\numberwithin{equation}{section}

\DeclareMathOperator{\ad}{Ad}

\DeclareMathOperator*{\spn}{span}
\DeclareMathOperator*{\clspn}{\overline{\spn}}

\newcommand{\cc}[1]{\mathcal{#1}}
\newcommand{\bb}[1]{\mathbb{#1}}

\newcommand{\N}{\bb N}
\newcommand{\Z}{\bb Z}
\newcommand{\Q}{\bb Q}

\newcommand{\T}{\bb T}
\renewcommand{\H}{\cc H}

\newcommand{\Chi}{\raisebox{2pt}{\ensuremath{\chi}}}

\newcommand{\inv}{^{-1}}

\newcommand{\eg}{\emph{e.g.}}
\newcommand{\cf}{\emph{cf.}}
\newcommand{\Ind}{\operatorname{Ind}}

\begin{document}

\title{Generalised Hecke algebras and $C^*$-completions}
\author[M. B. Landstad]{Magnus B. Landstad}
\address{Department of Mathematical Sciences, Norwegian University of Science
and Technology, N-7491, Trondheim, Norway.}
\email{magnusla@math.ntnu.no}
\author[N. S. Larsen]{Nadia S. Larsen}
\address{Department of Mathematics, University of Oslo,
P.O. Box 1053 Blindern, N-0316 Oslo, Norway.}
\email{nadiasl@math.uio.no} 

\thanks{This research was supported by the Research Council of Norway.}
\begin{abstract}
For a Hecke pair $(G, H)$ and a finite-dimensional representation
$\sigma$ of $H$ on $V_\sigma$ with finite range    we consider a generalised Hecke
algebra $\H_\sigma(G, H)$, which we study by embedding the given
Hecke pair in a Schlichting completion $(G_\sigma, H_\sigma)$ that
comes equipped with a continuous extension $\sigma$ of $H_\sigma$.
There is a (non-full) projection $p_\sigma\in C_c(G_\sigma, {\cc
B}(V_\sigma))$ such that $\H_\sigma(G, H)$ is isomorphic to
$p_\sigma C_c(G_\sigma, {\cc B}(V_\sigma))p_\sigma$. We study the
structure and properties of $C^*$-completions of the generalised
Hecke algebra arising from this corner realisation, and via
Morita-Fell-Rieffel equivalence we identify, in some cases explicitly,
the resulting proper ideals of $C^*(G_\sigma, {\cc B}(V_\sigma))$.
By letting $\sigma$ vary, we can compare these ideals.  The main
focus is on the case with $\dim\sigma=1$ and applications include
$ax+b$-groups and the Heisenberg group.
\end{abstract}

\maketitle

\section*{Introduction}

A Hecke pair $(G, H)$ consists of a group $G$ and a subgroup
$H$ such that $L(x):=[H:H\cap xHx^{-1}]$ is finite for all $x\in G$. Our interest lies
in studying $C^*$-completions of a \emph{generalised Hecke algebra} $\H_\sigma(G, H)$
associated with a Hecke pair $(G, H)$ and a  unitary representation $\sigma$ of $H$ on
a finite-dimensional Hilbert space $V_\sigma$. As a vector space,  $\H_\sigma(G, H)$ consists
of functions $f:G\to {\cc B}(V_\sigma)$ with
finite support in $H\backslash G/H$ such that
$$
f(hxk)=\sigma(h)f(x)\sigma(k) \text{ for all } h,k\in H, x\in G.
$$
When the group is locally compact totally disconnected, and the subgroup is compact and open, such
algebras, endowed with a natural convolution,
play a fundamental role in the representation theory
of reductive $p$-adic groups, see for example \cite{Ho}.

When $\sigma$ is the trivial representation of $H$, $\H_\sigma(G,
H)$ is the \emph{Hecke algebra}  $\H(G, H)$ of the pair $(G, H)$,
see e.g. \cite{K}. With an appropriate involution, $\H_\sigma(G, H)$
becomes a $*$-algebra. Our goal is to shed light on the structure of
$C^*$-completions of $\H_\sigma(G, H)$, and to identify conditions
which ensure that a largest $C^*$-completion exists. That this last
issue is important and non-trivial was demonstrated by Hall, who in
\cite{H} gave an example of a Hecke pair $(G, H)$ such that $\H(G,
H)$ does not have a largest $C^*$-completion.
We find it natural to
investigate the structure of $C^*$-completions in the more general
context of a generalised Hecke algebra $\mathcal{H}_\sigma(G, H)$ of
a Hecke triple $(G, H, \sigma)$. Results from \cite{KLQ} valid for
Hecke pairs $(G, H)$ are not directly applicable in the setup of
$(G, H,\sigma)$ with non-trivial $\sigma$, but the overall strategy
from \cite{KLQ} can be adapted and developed, as we shall show, in
order to deal with the differences that arise in our more general
context.

The interesting structure and
properties of the Hecke $C^*$-algebra
introduced  by Bost and Connes in \cite{BC}
have motivated intense research devoted to the
study of Hecke $C^*$-algebras of large classes of Hecke pairs, see for example
\cite{ALR, B,   BLPR, GW, H,  KLQ, LaLa0, LaR,  LarR, T}.
A powerful tool to analyse $\H(G, H)$ is the ``Schlichting completion''
$(\overline{G}, \overline{H})$ of $(G, H)$: this is a new Hecke pair
consisting of a locally compact totally disconnected group and a compact
open subgroup \cite{T}. Then $\H(G, H)$ is isomorphic to
$C_c(\overline{H}\backslash \overline{G}/\overline{H})$, which is a corner of
the group algebra $C_c(\overline{G})$, and this viewpoint
facilitates the analysis of $C^*$-completions in the realm of Banach $*$-algebras.
Tzanev's construction of $(\overline{G}, \overline{H})$ was inspired by work
of Schlichting, see for example \cite{Sch}, and was reviewed in \cite{KLQ},
where it is employed to study $C^*$-completions by looking at ideals in
$C^*(\overline{G})$,  see  \cite{GW, Laca_dil,  LarR} for other approaches.

In \cite{Cu0}, Curtis considers a Hecke algebra of a triple
$(G, H, \sigma)$ where $\sigma$ is a finite-dimensional
unitary representation of $H$, and studies a von Neumann algebra naturally
associated to it. Curtis constructs a completion
$(G_\sigma, H_\sigma,\sigma)$, but does not prove that it is unique.

In the present study we concentrate the attention to the case where the
representation $\sigma$ of $H$ has  finite range.
 There are several reasons for this: we
require finite range first because this case is simpler to handle, but more importantly,
the case with $\sigma(H)$ infinite is fundamentally different, since the completion of
$G$ then will contain a copy of $\T$ and therefore will not be totally disconnected.

With our extended theory  we use both Fell's and Rieffel's versions
of Morita equivalence to analyse the structure of $C^*$-completions
of $\H_\sigma(G, H)$ with respect to $\sigma$. It turns out that we
compare corners of $C_c(G_\sigma, {\cc B}(V_\sigma))$ determined by projections $p_\sigma$.
If $\sigma$ is the trivial representation, it was shown in
\cite{KLQ} that the projection $p_\sigma$ often is full. However, for nontrivial $\sigma$ it turns out that $p_\sigma$ is never full.

In the by now
classical example of the Bost-Connes Hecke pair, the completion
$G_\sigma$ will be the same for all finite characters $\sigma$ of
$H$, and therefore all the projections $p_\sigma$ live in one group
$C^*$-algebra $A=C^*(G_0)$, see Example~\ref{ex:full ax+b}. It follows that the generalised Hecke
$C^*$-algebras $p_\sigma A p_\sigma$ are all Morita-Rieffel
equivalent to ideals in the same $C^*$-algebra $A$, and can
therefore be compared more naturally.  It turns out that these
ideals are built from the primitive ideals of the Bost-Connes Hecke
$C^*$-algebra identified by Laca and Raeburn in \cite{LaR2}.

The organisation of the paper is as follows. In section~\ref{start} we construct
our Schlichting completion  $(G_\sigma,
H_\sigma, \sigma)$ of a Hecke triple $(G, H,\sigma)$, where $\sigma$ is a
finite-dimensional unitary representation of $H$ with finite range. We prove in
Theorem~\ref{uniqueness_univ_prop} that $(G_\sigma, H_\sigma, \sigma)$
has a universal property, which is essentially provided
by the universal property of Schlichting completions of Hecke pairs, see
\cite[Theorem 3.8]{KLQ}. In section~\ref{section_gen_Hecke} we define the
generalised Hecke algebra $\H_\sigma(G, H)$ for arbitrary finite-dimensional
$\sigma$ and we realise $\H_\sigma(G, H)$ as a corner of
$C_c(G_\sigma, {\cc B}(V_\sigma))$.
 One  new ingredient that appears in the study of generalised as opposed
to usual Hecke algebras is that not every double coset $HxH$ for
$x$ in $G$ supports a non-zero function in $\H_\sigma(G, H)$. When $\sigma$ is
one-dimensional, we obtain further insight. The subset $B$ of
the $x$'s in $G$ which do support a non-zero function need not be a subgroup of
$G$, but nevertheless it harmonises with the Schlichting completion. If $B$ is a group,
its closure $B_\sigma$ in $G_\sigma$ determines a corner $p_\sigma C_c(B_\sigma)p_\sigma$,
and we prove (Proposition~\ref{H_sigma_spanned_by_epsilonx})
that this corner is isomorphic to $\H_\sigma(G, H)$.

Section~\ref{vN_alg} contains an analysis of the continuity properties of the
induced representation from $H$ to $G$ with respect to the process of
taking the Schlichting completion. To illustrate the point that our
Schlichting completion of $(G, H, \sigma)$ is a
profitable alternative to studying $\H_\sigma(G, H)$, we employ it to give a
short proof (see Theorem~\ref{commutant_is_Hecke_vNalg}) of a classical result
which asserts that the commutant of the
induced representation $\Ind_H^G\sigma(G)$ is the weak closure of the
``intertwining operators'', cf. \cite[Theorem 2.2]{Bi1} (which mends the
apparently deficient proof of \cite[Theorem 3]{Co}) or \cite[Proposition 1.3.10]{Cu1}.
As an immediate corollary, for one-dimensional $\sigma$
we describe the irreducibility of $\Ind_H^G\sigma$ in terms of the Hecke algebra, thus
recovering Mackey's condition in \cite{M}.

We shall  often specialise to  one-dimensional representations, and in section~\ref{completions} we start by showing
 that if $\dim \sigma>1$, the
generalised Hecke algebra is still Morita equivalent to the ideal in
$C_c(G_\sigma)$ generated by the character of $\sigma$.
 So also in this case the generalised Hecke algebras can be studied by
looking at ideals in $C^*(G_\sigma)$.
For these ideals we describe the
nondegenerate representations as in \cite[\S 5]{KLQ} by means of a category equivalence,
see Corollary~\ref{category_equiv}.

In the presence of a normal subgroup $N$ of
$G$ which contains $H$ we describe the structure of these ideals as twisted
crossed products.
 If in addition $H$ is normal in $N$ and $N\subset
B$, we can conclude that $\H_\sigma(G, H)$ has a largest
$C^*$-completion, see
Corollary~\ref{C*_completion_normal_subgroups}. Finally, we study the special, but
interesting instance where $B$ is a group and $(B, H)$ is directed
in the sense of \cite[\S 5]{KLQ}. A largest $C^*$-completion turns
out to exist, and we can give a concrete description of the ideal in
$C^*(G_\sigma)$, see Theorem~\ref{theo_p1_B_directed}.

The last section is devoted to applications.
We show in Proposition~\ref{sigma_extends} that $\H_\sigma(G, H)$ is isomorphic
to $\H(G, H)$ when $\sigma$ extends to a character of $G$.
 We illustrate in examples the multitude of possible outcomes of the construction of the
Schlichting completion.

\smallskip
\noindent \textbf{Conventions.}
For a  Hecke pair $(G, H)$ and a subset $X$ of $G$, the
notation $y\in X/H$ (and $y
\in H\backslash{X}/H$)  means that $y$ runs over a set of  representatives for
the left cosets $X/H$ (and the double cosets $H\backslash{X}/H$).

All representations of topological groups are assumed to be unitary and continuous.
If $L$ is a locally compact group with a left invariant Haar measure $\mu$
and modular function $\Delta$, then the space $C_c(L)$ of compactly supported
continuous functions on $L$ is a $*$-algebra with multiplication given by
usual convolution $f\ast g(x)=\int_L f(y)g(y^{-1}x)d\mu(y)$,
and involution given by $f^*(x)=\Delta(x^{-1})\overline{f(x^{-1})}$.
The  group $C^*$-algebra $C^*(L)$ is generated by a universal unitary
representation of $L$ into the unitary group of the multiplier algebra
$M(C^*(L))$, and $\{\int_L f(x)x d\mu(x)\mid f\in C_c(L)\}$
spans a dense subspace of $C^*(L)$, where we identify $x$ in $L$ with its image in
$M(C^*(L))$.
\medskip

\section{Hecke pairs and Schlichting completions}\label{start}

We recall from \cite[Definition 3.3]{KLQ} that
if $(G, H)$ is a Hecke pair, then the collection $\{xHx^{-1}\mid x\in G\}$ is a
neighbourhood  subbase for the Hecke topology on $G$ from $(G, H)$.
If a Hecke pair $(G, H)$ is such that
\begin{equation}
\bigcap_{x\in G}xHx^{-1}=\{e\},
\label{reduced_pair}
\end{equation}
then it is called \emph{reduced} \cite{T}, and the Hecke topology from $(G, H)$ is
Hausdorff. The \emph{Schlichting completion} of a Hecke pair
$(G, H)$ was constructed in \cite{T} to be an essentially unique Hecke pair consisting of a
locally compact totally disconnected group with a compact open subgroup in which $G$ and
respectively $H$ embedd densely.
In the terminology of  \cite[\S 3]{KLQ}, the Schlichting completion
of $(G, H)$ consists of the closures of $G$ and $H$ in
the Hecke topology from $(G, H)$. The Schlichting completion of a Hecke pair is a
\emph{Schlichting pair}, which by \cite[\S 3]{KLQ} is a
reduced Hecke pair with the additional feature that
the underlying subgroup is compact and open in its corresponding Hecke
topology.

Suppose that $(G, H)$ is a Hecke pair and $\sigma$ is a finite-dimensional unitary
representation of $H$ on a Hilbert space $V_\sigma$ such that $\sigma(H)$ is finite.
Then $K:=\ker \sigma$ is a normal subgroup of $H$ of finite index,
and hence $(G, K)$ is a Hecke pair. Let $(G_\sigma, K_\sigma)$ denote the Schlichting
completion of $(G, K)$. We have the following lemma.

\begin{lem}
\textnormal{(a)}
 The closure $H_\sigma$ of $H$ in the Hecke topology from $(G, K)$ is a compact open subgroup
 of $G_\sigma$.

\smallskip
\textnormal{(b)}
 $\sigma$ is continuous for the Hecke topology from  $(G, K)$,
and thus has a unique extension to a finite-dimensional unitary representation $\sigma$ of
$H_\sigma$ with kernel $K_\sigma.$
\label{construction_of_Hecke_triple}
\end{lem}

\begin{proof}
We claim that $hK\to h K_\sigma$ for $h\in H$ is an isomorphism
$
H/K\overset{\cong}{\longrightarrow} H_\sigma/K_\sigma;
$
indeed, an element $h\in H\setminus K$ is carried to the
open set $hK_\sigma$ which is disjoint from $K$,
showing injectivity, and
surjectivity follows because any given $xK_\sigma$ in
$H_\sigma/K_\sigma$ is open, and hence meets the dense subset $H$ of
$H_\sigma$.  Thus $H_\sigma/K_\sigma$ is finite, so
$H_\sigma$ is compact because its quotient by a compact subgroup is
again compact. This proves (a). For (b) it suffices to show
continuity of $\sigma$ at $e$, and this follows by inspection using that
$\sigma(H)$ is finite and $K$ is open.
\end{proof}

\begin{defn} A \emph{$n$-dimensional Hecke triple} $(G, H, \sigma)$ consists of a Hecke pair
$(G, H)$ and a $n$-dimensional unitary representation $\sigma$ of $H$ on $V_\sigma$
with finite range. We say  that $(G, H, \sigma)$
is \emph{reduced} if the Hecke pair $(G, K:=\ker \sigma)$
is reduced. We call the Hecke triple $({G}_\sigma, {H}_\sigma, \sigma)$ from
Lemma~\ref{construction_of_Hecke_triple}
the \emph{Schlichting completion} of $(G, H, \sigma)$.
\label{def_of_Hecke_top_from_sigma}
\end{defn}

\begin{rem} Suppose that $L$ is a totally
disconnected locally compact group, $M$ a compact open subgroup, and
$\rho$ a continuous finite-dimensional unitary representation of
$M$. It follows from \cite[Corollary (28.19)]{HR} that $\rho(M)$ is
finite.
\end{rem}

We next prove that the Schlichting completion of a Hecke triple
has a universal property, and is unique up to topological isomorphism.

\begin{thm}\label{uniqueness_univ_prop}
Let $(G, H, \sigma)$ be a reduced $n$-dimensional  Hecke triple.
Then the Schlichting completion $({G}_\sigma, {H}_\sigma, {\sigma})$
has the following universal property: suppose that $L$ is a locally
compact totally disconnected group, $M$ a compact open subgroup,
$\rho$ a unitary representation of $M$ on $V_\sigma$, and $\phi:G\to
L$ a homomorphism such that $(L, \ker \rho)$ is reduced,
 $\phi(G)$ is dense in $L$, $\phi(H)\subseteq M$
and $\sigma=\rho \circ \phi\vert_H$. Then there is a unique
continuous homomorphism
$\overline{\phi}$ from ${G}_\sigma$ onto $L$ which extends $\phi$ and satisfies the
identity
\begin{equation}
\rho \circ \overline{\phi}\vert_{{H}_\sigma}={\sigma}.
\label{phi_extended}
\end{equation}

If in addition $\phi^{-1}(M)=H$, then $\overline{\phi}$ will be a topological group
isomorphism of ${G}_\sigma$ onto $L$ and of ${H}_\sigma$ onto $M$.
\end{thm}

\begin{proof}
Denote $N:=\ker \rho$. Then $\phi(K)\subseteq N$.
Applying the first half of  \cite[Theorem 3.8]{KLQ} to the Schlichting pair
$(L, N)$
gives a unique continuous homomorphism $\overline{\phi}$
from $G_\sigma$ into $L$ which extends $\phi$. Since
$$
\rho \circ \overline{\phi}\vert_H(h)=\rho\circ \phi(h)=\sigma(h)
$$
for all $h\in H$, \eqref{phi_extended} follows from the continuity of $\rho \circ \overline{\phi}$ and $\sigma$ on $H_\sigma$.

If $\phi^{-1}(M)=H$, a straightforward verification then shows that $\phi^{-1}(N)
=K$, and it follows from the second half of \cite[Theorem 3.8]{KLQ} that
$\overline{\phi}$ is a topological group isomorphism of $G_\sigma$ onto $L$. The assumption
$\phi^{-1}(M)=H$ implies that $\phi(H)=M\cap \phi(G)$. Since $M$ is open and closed,
\begin{equation}
\overline{\phi(H)}=\overline{M\cap \phi(G)}=M\cap \overline{\phi(G)}=M.\label{M_is_phiH_closed}
\end{equation}
The set $\overline{\phi}(H_\sigma)$ is compact, hence closed, and so it equals
$\overline{\phi(H)}$. By invoking  \eqref {M_is_phiH_closed} we obtain the last claim of the theorem.
\end{proof}

\begin{rem}\label{def_of_iota}
Let $(G, H, \sigma)$ be a reduced Hecke triple, $(G_\sigma, H_\sigma, \sigma)$
the Schlichting completion, and $j_1$ the dense embedding of $G$ in
$G_\sigma$. Denote by $(G_0, H_0)$
the Schlichting completion of $(G, H)$ and  by $j_0$ the dense embedding
$G\to G_0$. Since
$K\subseteq H\subseteq H_0$, the first half of
\cite[Theorem 3.8]{KLQ} gives
a continuous homomorphism $\iota:{G}_\sigma \to G_0$ such that $\iota\circ j_1=j_0$.
Since $H$ is dense in both $H_\sigma$ and $H_0$, $\iota(H_\sigma)=H_0$.
We typically omit $j_0$ and
$j_1$ from the notation.
\end{rem}

\begin{rem} In \cite{Cu0}, for a Hecke pair $(G, H)$ and a finite-dimensional
unitary representation $\sigma$ of $H$, Curtis defines
an equivalence relation $\sim$ on $G\times \sigma(H)$  by
$(g, t)\sim (gh^{-1}, \sigma(h)t)$ for all $h\in H$. With
$S_\sigma:=(G\times \sigma(H))/\sim$ denoting the quotient space,
$G$ is endowed with the topology pulled back from the
compact-open topology on the space of
continuous functions $\{f:S_\sigma\to S_\sigma\}$. Then part of
\cite[Theorem 3]{Cu0} asserts that the closures of $G$ and $H$ in this topology
and the unique extension
of $\sigma$ to the closure of $H$ have the universal property. If
$\sigma$ has finite range, note that the map $g\mapsto [g, 1]$ from $G$ onto
$S_\sigma$  is a bijection from $G/K$
onto $S_\sigma$, which is equivariant for the actions of $G$ as permutations on
$G/K$ and on $S_\sigma$. Thus $G_\sigma$ is the same as the completion
constructed in \cite{Cu0}.
\end{rem}

\section{The generalised Hecke algebra of $(G, H, \sigma)$}\label{section_gen_Hecke}

The next definition appears in \cite{Cu0}, with the difference that $H$ is an
arbitrary subgroup of $G$ and one takes functions $f$ with  finite support on
$H\backslash G$ and $G/H$.  However, for  the purposes
of using Schlichting completions, the important case, also in
\cite{Cu0}, is that of a Hecke pair $(G, H)$.

\begin{defn}
Given a $n$-dimensional Hecke triple  $(G, H,\sigma)$,
 let $\H_\sigma(G, H)$ be the vector space of functions $f:G\to {\cc B}(V_\sigma)$ which  have
finite support in $H\backslash G/H$ and satisfy $f(hxk)=\sigma(h)f(x)\sigma(k)$
for all $h,k\in H, x\in G$. The
\emph{generalised Hecke algebra} associated with $(G, H,\sigma)$ is
$\H_\sigma(G, H)$ endowed with the convolution
\begin{equation}
f\ast g(x)=\sum_{yH\in G/H}f(y)g(y^{-1}x).\label{def_of_conv}
\end{equation}
The identity element is the function $\varepsilon_H$ defined by $\varepsilon_H(x)=
\sigma(x)$ when $x\in H$ and $\varepsilon_H(x)=0$ otherwise.
\end{defn}

The key reason that motivates the study of $\H_\sigma(G, H)$
in terms of the Schlichting completion $(G_\sigma, H_\sigma, \sigma)$ of $(G, H,
\sigma)$ is the following standard result.

\begin{lem}\label{gen_Hecke_on_lcg}
Let $L$ be a locally compact group, $M$ a compact open subgroup,
$\rho$ a  finite-dimensional
unitary representation  of $M$,
and choose the Haar measure $\mu$ on $L$ normalised so that $\mu(M)=1$.
Then $\H_\rho(L, M)$ is equal to the subalgebra
\begin{equation}
\{f\in C_c(L,{\cc B}(V_\rho))\mid f(mxn)=\rho(m)f(x)\rho(n), \forall m,n\in M, x\in L\}\label{Hecke_algebra_as_functions_in_Cc}
\end{equation}
of $C_c(L,{\cc B}(V_\rho))$, endowed with the convolution with respect to $\mu$.
\end{lem}

\begin{prop} Suppose that $(G, H, \sigma)$ is a reduced Hecke triple. Let
$(G_\sigma, H_\sigma, \sigma)$ be the Schlichting completion of $(G, H, \sigma)$, and choose
the Haar measure $\mu$ on $G_\sigma$ normalised so that $\mu(H_\sigma)=1$. Then the map
$\Psi:\H_\sigma({G_\sigma}, {H_\sigma})\to \H_\sigma(G,
H)$ given by $\Psi(f)=f\vert_G$  is an algebra isomorphism.
\label{Hecke_algebra_for_G_and_Gbar_same}
\end{prop}

\begin{proof}
By  adapting the argument in \cite[Proposition 3.9 (iii)]{KLQ} to the reduced pair $(G, K)$, it follows that
$HxH\mapsto H_\sigma xH_\sigma$ for $x\in G$ is a bijection from
$H\backslash G/H$ onto $H_\sigma\backslash G_\sigma /H_\sigma$. Since $G_\sigma$
and $H_\sigma$ contain dense copies of $G$ and $H$, it follows that
the map $\Psi$ is well-defined.

Given $f$ in $\H_\sigma(G, H)$, note that by the invariance property of $f$,
$$
f((xKx^{-1})x)=f(xK)=f(x)\sigma(K)=f(x)
$$
for all $x\in G$. Thus $f$ is continuous for the Hecke topology from $(G, K)$,
and so extends to a function in $\H_\sigma (G_\sigma, H_\sigma)$. It follows that
 $\Psi$ is bijective. A routine calculation shows that $\Psi(f\ast g)=\Psi(f)\ast
\Psi(g)$, and the claim follows.
\end{proof}

From Proposition~\ref{Hecke_algebra_for_G_and_Gbar_same} and
Lemma~\ref{gen_Hecke_on_lcg} it seems natural to define the
involution on $\H_\sigma(G, H)$  as in $C_c(G_\sigma, {\cc B}(V_\sigma))$
by using the modular function of $G_\sigma$, so we investigate
its meaning for the original Hecke triple.
This is in fact answered by Schlichting in \cite[Lemma 1(iii)]{Sch}. We
need some notation first. For a Hecke pair $(G, H)$ and any $x$ in $G$ let
$H_x:=H\cap xHx^{-1}$, $L(x):=[H:H_x]$ and $\Delta_H(x):= {L(x)}/{L(x^{-1})}$.
By \cite{Sch},  if $H$ is a compact open subgroup of a locally compact group
$G$ the modular function $\Delta$ of $G$ satisfies $\Delta(x)=
\Delta_H(x)$. In particular, $\Delta_H(x)$ does not depend on which compact
open subgroup we use. With the notation of Remark~\ref{def_of_iota} we obtain:

\begin{cor}\label{modular_functions_same}
The modular functions $\Delta_\sigma$ of $G_\sigma$ and
$\Delta_0$ of $G_0$  satisfy $\Delta_0\circ\iota=\Delta_\sigma$.
\end{cor}

If  $(G, H, \sigma)$ is a  reduced Hecke triple with Schlichting completion
$(G_\sigma, H_\sigma, \sigma)$, then
$[H:H_x]=[H_\sigma:(H_\sigma)_x]$, so we can supplement
Proposition~\ref{Hecke_algebra_for_G_and_Gbar_same}.

\begin{prop}\label{Hecke*algebra_for_G_and_Gbar_same}
Let  $(G, H,\sigma)$ be a reduced Hecke triple with $K:=\ker \sigma$, and
define an involution on
$\H_\sigma(G, H)$ by
\begin{equation}
f^*(x)=\Delta_K(x^{-1})f(x^{-1})^*, \text{ for }x\in G.
\label{def_involution}
\end{equation}
Then the map
$\Psi$ of Proposition~\ref{Hecke_algebra_for_G_and_Gbar_same}
is an isomorphism of $*$-algebras.
\end{prop}

\begin{rem}  Most authors do not include $\Delta$ in the definition of an
involution on a (generalised) Hecke algebra.
But we claim that this is more natural,
for instance the $l^1$-norm on $\H_\sigma(G, H)$ defined by
\begin{equation}
\Vert f\Vert_1=\sum_{y\in G/H}\| f(y)\| \text{ for }f\in \H_\sigma(G, H)
\label{def_of_norm}
\end{equation}
satisfies $\| f^*\|_1=\| f\|_1$. We let
$l^1(G, H, \sigma)$ be the completion  of $\H_\sigma(G, H)$ in $\Vert \cdot
\Vert_1$. As a consequence of Proposition~\ref{Hecke_algebra_for_G_and_Gbar_same}
and Proposition~\ref{Hecke*algebra_for_G_and_Gbar_same} we get the following:
\end{rem}

\begin{prop} With the assumptions and notation from
Proposition~\ref{Hecke_algebra_for_G_and_Gbar_same}, the map $\Psi$
extends to an isomorphism $L^1({G_\sigma}, {H_\sigma},
{\sigma})\cong l^1(G, H, \sigma)$ of Banach $*$-algebras.
\label{L1_algebras_isom}
\end{prop}

\begin{proof}
A computation shows that $\Vert \Psi(f)\Vert_1=\Vert f\Vert_1$ for all $f\in
\H_\sigma({G_\sigma}, {H_\sigma})$. Hence $\Psi$ extends
to the completion in the norm from (\ref{def_of_norm}).

Note that (\ref{def_of_norm}) on $\H_\sigma({G_\sigma},
{H_\sigma})$ is the usual $L^1$-norm on
$C_c({G_\sigma, {\cc B}(V_\sigma)})$. A routine
calculation shows that $L^1({G_\sigma},
{H_\sigma}, {\sigma})$ is a closed $*$-subalgebra of
$L^1({G_\sigma, {\cc B}(V_\sigma)})$ for the natural involution of $L^1({G_\sigma, {\cc B}(V_\sigma)})$.
Hence $L^1({G_\sigma}, {H_\sigma}, {\sigma})$ is a Banach $*$-algebra and the claim follows.
\end{proof}

\begin{thm} Given a reduced $n$-dimensional Hecke triple $(G, H, \sigma)$, let $(G_\sigma, H_\sigma,
\sigma)$ denote its Schlichting completion. Denote by $\mu$ the left invariant
Haar measure on ${G_\sigma}$ such that $\mu({H_\sigma})=1$.
Then the function $p_\sigma(x):=\Chi_{{H_\sigma}}(x)
{\sigma}(x)$ is a self-adjoint projection
in $C_c({G_\sigma, {\cc B}(V_\sigma)})$, and we have isomorphisms
between the $*$-algebras $ p_\sigma  C_c(G_\sigma, {\cc
B}(V_\sigma))  p_\sigma=\H_\sigma(G_\sigma, H_\sigma)$ and
$\H_\sigma(G, H)$, and between the  Banach $*$-algebras $p_\sigma
L^1(G_\sigma, {\cc B}(V_\sigma)) p_\sigma$ and $l^1(G, H, \sigma).$
\label{fundam_proj}
\end{thm}

\begin{proof}
It is straightforward that $p_\sigma$ is a self-adjoint  projection.
Lemma~\ref{gen_Hecke_on_lcg} implies that
$\H_\sigma(G_\sigma, H_\sigma)=p_\sigma C_c(G_\sigma, {\cc B}(V_\sigma))p_\sigma$, and then the
claimed isomorphisms follow from Proposition~\ref{Hecke*algebra_for_G_and_Gbar_same}
and Proposition~\ref{L1_algebras_isom}.
\end{proof}

\section{The generalised Hecke algebra of $(G, H, \sigma)$ when $\dim\sigma=1$}
\label{section_gen_Hecke_dim=1}

In this section we assume $\dim(V_\sigma)=1$, although some of the results still will be true
in greater generality, in particular the results about the set $B$, \cf\  \cite{Bi1, Bi2}.

Similarly to \cite[Lemma 4.2(iii)]{KLQ} we have that $p_\sigma C_c(G_\sigma)
p_\sigma=\spn\{p_\sigma x p_\sigma\mid x\in G_\sigma\}$. We proceed to
identify functions in a spanning set for
$\H_\sigma(G, H)$ which correspond by Theorem~\ref{fundam_proj} to the products
$p_\sigma xp_\sigma$.

It is known, see for instance \cite{K}, that the Hecke algebra
of a pair $(G, H)$ is linearly spanned by the collection $\{\Chi_{HxH}\mid
x\in H\backslash G/H\}$ of characteristic functions of double cosets.
To account for non-zero functions in $\H_\sigma(G, H)$ supported on a
given double coset, note that for $f\in \H_\sigma(G, H)$, $x\in G$ and
$h\in H_x=H\cap xHx^{-1}$ we have
$$
\sigma(h)f(x)=f(hx)=f(xx^{-1}hx)=f(x)\sigma(x^{-1}hx).
$$
Thus for $f$ to be supported on $HxH$ we need $
\sigma(h)=\sigma(x^{-1} hx)\text{ for }h\in H_x$. This is condition $(t_g)$ in
\cite[Proposition 1.2]{Bi1}, and goes at least back to Mackey \cite{M}. We denote
\begin{equation}
B:=\{x\in G\mid \sigma(h)=\sigma(x^{-1} hx)\text{ for }h\in H_x\}.
\label{special_x}
\end{equation}

\begin{rem}
The set $B$ contains $H$, is  closed under inverses,
and satisfies $BH=B=HB$. In general, $B$ is not a group, see for example
\cite[Example 1.4.4]{Cu1}. However, $B$ is a group in many cases such as, for
instance, when $\sigma$ extends to a character of $G$.
\end{rem}

\begin{lem}\label{basis_for_Hecke_algebra} For each $x\in B$ there is a well-defined element
\begin{equation}
{\varepsilon_x}(y)=
\begin{cases}
\sigma(hk)&\text{ if }y\in HxH\text{ and }y=hxk\\
0&\text{ if }y\notin HxH
\end{cases}
\label{def_of_epsilon_x}
\end{equation}
in $\H_\sigma(G, H)$, and the set $\{\varepsilon_x\mid x\in
H\backslash {B}/H\}$ forms a linear basis for $\H_\sigma(G, H)$.
\end{lem}

\begin{proof}
Clearly $\varepsilon_x\in \H_\sigma(G, H)$ and is well-defined
(these functions are essentially the elementary intertwining
operators from \cite{Bi1}). Since $\varepsilon_{h_0xk_0}(hxk)=
\overline{\sigma(h_0)\sigma(k_0)}\varepsilon_x(hxk)$ for all $h_0,
k_0\in H$,  different choices of representatives for the double
coset $HxH$ give rise to functions $\varepsilon_{h_0xk_0}$ which are
scalar multiples of $\varepsilon_x$, and since distinct double
cosets do not support a common $\varepsilon_x$ the lemma follows.
\end{proof}

When $B$ is a group, $(B, H,\sigma)$ is a new Hecke triple, we trivially have
$\H_\sigma(G, H)=\H_\sigma(B, H)$, and $\H_\sigma(G, H)$ has a linear basis
indexed over the double cosets of $B$ with respect to its subgroup $H$. To study $\H_\sigma(G, H)$
in this case we must naturally view it as a corner in $C_c(B_0)$ with $B_0$ denoting
the Schlichting completion of $(B, K)$. As  example~\ref{ex:full ax+b} shows, the Schlichting completion of
$(B, H,\sigma)$ need not come from $(G_\sigma, K_\sigma)$, because on $B$ the Hecke topology from $(G, K)$
differs from the Hecke topology from $(B, K)$. Nevertheless, the corner in $C_c(B_0)$
which is isomorphic to $\H_\sigma(B, H)$ is completely determined by the topology on $G_\sigma$. To prove
this, we first establish that the closure of
$B$ in $G_\sigma$ is precisely the set
defined by \eqref{special_x} for $(G_\sigma, K_\sigma)$.

\begin{lem}\label{closure_of_Hx}
We have $\overline{H\cap xHx^{-1}}=H_\sigma \cap xH_\sigma x^{-1}$ in $G_\sigma$ for all $x\in G$.
\end{lem}

\begin{proof}
Since $H_x$ is included in the closed set $(H_\sigma)_x$ for every $x\in G$,
we obtain one inclusion. Suppose that $h\in (H_\sigma)_x$. Let $F$ be a finite
subset of $G$, and take
$K_{\sigma, F}=\bigcap_{y\in F} yK_\sigma y^{-1}$, a neighbourhood of $e$. By
restricting, if needed, to a smaller neighbourhood, we may assume that $x\in F$.
Since $H$ is dense in $H_\sigma$, it intersects the open neighbourhood
$hK_{\sigma, F}$ of $h$. Thus we
have
\begin{equation}
hk_1=h_1, hk_2=xh_2x^{-1},\label{many_h}
\end{equation}
with $h_1, h_2\in H$ and $k_1, k_2\in K_{\sigma, F}$. Then $k_1^{-1}k_2=
h_1^{-1}xh_2x^{-1}$ is an element of $G\cap j_1^{-1}(K_{\sigma, F})$, which is
$\bigcap_{y\in F} yKy^{-1}$ because the Schlichting completion satisfies $j_1^{-1}
(K_\sigma)=K$. Thus $k_1^{-1}k_2$ lies in $xHx^{-1}$, and so
$$
h_1=xh_2x^{-1}k_2^{-1}k_1\in xHx^{-1}\cap H,
$$
from which it follows via \eqref{many_h} that $h_1\in hK_{\sigma, F}\cap H_x$.
Since this holds for all neighbourhoods $K_{\sigma, F}$, we have
$h\in \overline{H_x}$, as claimed.
\end{proof}

\begin{prop} Let $(G, H, \sigma)$ be a reduced $1$-dimensional Hecke triple, and consider its
Schlichting completion $(G_\sigma, H_\sigma, \sigma)$. Let
$B_\sigma$ denote
\begin{equation}
\{x\in  G_\sigma\mid \sigma(h)=\sigma(x^{-1}hx)\text{ for }h\in H_\sigma \cap
xH_\sigma x^{-1}\}.\label{def_of_Bsigma}
\end{equation}
Then $B_\sigma$ is equal to the closure of $B$ in the Hecke topology from $(G, K)$.
\label{Hx_Kx}
\end{prop}

\begin{proof}
If $(x_i)$ is a net in $B_\sigma$ converging to $x$, then eventually $H_\sigma \cap x_iH_\sigma x_i^{-1}$
coincides with $H_\sigma \cap xH_\sigma x^{-1}$, and so $B_\sigma$
is closed. Since $B_\sigma H_\sigma=B_\sigma$, it is also open.

Lemma~\ref{closure_of_Hx} implies that
$$
B_\sigma\cap G=\{x\in  G\mid \sigma(h)=\sigma(x^{-1}hx)\text{ for }h\in
H_\sigma \cap xH_\sigma x^{-1}\}=B.
$$
Hence the closure of $B$ is included in $B_\sigma$. To show
equality, take $x\in B_\sigma$ and $K_{\sigma, F}$  a neighbourhood of $e$. We
must show that $xK_{\sigma, F}$ has non-empty intersection with $B$. By density
of $G$ in $G_\sigma$ there is $k\in K_{\sigma, F}$ such that $xk\in G$. From
$K_{\sigma, F}\subset H_\sigma$ and $B_\sigma H_\sigma =B_\sigma$ it follows that
$xk\in B_\sigma\cap G=B$, as claimed.
\end{proof}

\begin{prop}\label{H_sigma_spanned_by_epsilonx}
 Let $(G, H, \sigma)$ be a reduced $1$-dimensional Hecke triple, let
$(G_\sigma, H_\sigma, \sigma)$ be its Schlichting completion, and
assume that $B$ is a subgroup of $G$. Let $B_\sigma$ be the subgroup
of $G_\sigma$ defined in \eqref{def_of_Bsigma}. Then $\H_\sigma(G,
H)$ is isomorphic to $p_\sigma C_c(B_\sigma)p_\sigma$.
\end{prop}

\begin{proof}
Since $B_\sigma$ is closed in $G_\sigma$, it is locally compact.
The subgroup $H_\sigma$ is open and compact in $B_\sigma$, so on one hand
$B_\sigma$ and $G_\sigma$ have
the same modular function, equal to $\Delta_{K_\sigma}$ and $\Delta_{H_\sigma}$,
and on the other  $\H_\sigma(B_\sigma, H_\sigma)$  equals
$p_\sigma C_c(B_\sigma)p_\sigma$  by Lemma~\ref{gen_Hecke_on_lcg}.

Given $x$ in $B_\sigma$, it follows from Proposition~\ref{Hx_Kx} that we
can pick $b$ in $B$ such that $xH_\sigma=j_1(b)H_\sigma=\overline{j_1(bH)}$. Thus
the isomorphism  $\Psi$ from
Proposition~\ref{Hecke*algebra_for_G_and_Gbar_same} carries functions supported on
$H_\sigma xH_\sigma$ with $x\in B_\sigma$ to functions supported on $HbH$ with
$b\in B$. Hence by  Lemma~\ref{basis_for_Hecke_algebra} the map $\Psi$ is an
isomorphism of $\H_\sigma(B_\sigma, H_\sigma)$ onto $\H_\sigma(B, H)$.
\end{proof}

\begin{lem}
\label{psigma_x_psigma} With the notation of Theorem~\ref{fundam_proj} we have
\begin{equation}
p_\sigma xp_\sigma
=\begin{cases}\frac 1{L(x)}\varepsilon_x&\text{ if }x\in B_\sigma\\
0&\text{otherwise}.
\end{cases}
\label{p_sigma_corner_as_epsilon_x}
\end{equation}
In particular, $\Vert p_\sigma xp_\sigma\Vert_1=1$ for every $x$ in $B_\sigma$.
\end{lem}

\begin{proof}
Since $p_\sigma$ is supported on $H_\sigma$, the product $p_\sigma xp_\sigma$
is supported on the double coset $H_\sigma x H_\sigma$.  The claim then follows
because
\begin{align}
p_\sigma xp_\sigma(hxk)
&=\sigma(h)\sigma(k)\int_{H_\sigma\cap xH_\sigma x^{-1}}\sigma(l)\overline{\sigma(x^{-1}lx)}  dl \notag \\
&=\begin{cases}
\sigma(h)\sigma(k)\mu(H_\sigma\cap xH_\sigma x^{-1})&\text{ if }x\in B_\sigma\\
0&\text{ otherwise}.\notag\\
\end{cases}
\end{align}
\end{proof}

\section{The von Neumann algebra of $(G, H,\sigma)$ and induced representations}\label{vN_alg}

One motivation for studying (generalised) Hecke algebras is that
they usually generate the commutant of the corresponding induced
representation. In particular this gives irreducibility criteria for
induced representations of locally compact groups, see \eg\ Mackey
\cite{M}, Corwin \cite{Co} and Binder \cite{Bi1, Bi2}. In
\cite[Theorem 2.2]{Bi1}, the commutant of the induced representation
is realised  as the weak closure of the algebra of so-called
elementary intertwining operators (this result was obtained earlier
by Corwin \cite{Co}, but that proof is claimed to be incomplete
\cite{Bi1}), see also \cite[Proposition 1.3.10]{Cu1}. In this
section we study the continuity properties of the induced
representation arising from a Hecke triple, and  obtain as a direct
consequence a new proof of Binder's theorem for such
representations.

Suppose that $(G, H, \sigma)$ is a $n$-dimensional Hecke triple. The induced representation
$\lambda_\sigma:=\Ind_H^G \sigma$ acts by the formula
$\lambda_\sigma(x)(f)(y):=f(x^{-1}y)$
  in the space $l^2(G, H, \sigma)$ defined as
\begin{equation}
\{f:G\to V_\sigma\mid f(xh)={\sigma(h^{-1})}f(x),\forall x\in G, h
\in H, \text{ and }\sum_{x\in G/H} \| f(x)\|^2<\infty\}.
\label{Hilbert_space_of ind_rep}
\end{equation}
The function  $\delta_{y,\xi}:G\to V_\sigma$ defined by
$\delta_{y,\xi}(z)=\Chi_H(z^{-1}y) \sigma(z^{-1}y)\xi$ for $y\in G $
and $\xi\in V_\sigma$ lies in $l^2(G, H, \sigma)$. If we choose a
set of coset representatives $y\in G/H$ and an orthonormal basis
$\{\xi_i\mid i=1,\cdots, n\}$ for $V_\sigma$, then
$\{\delta_{y,\xi_i}\mid y\in G/H, i=1,\cdots, n\}$ is an orthonormal
basis for $l^2(G, H, \sigma)$.

\begin{lem}\label{induced_rep_extends} Suppose that $(G, H, \sigma)$ is a reduced
Hecke triple. Let $K=\ker \sigma$ and $({G}_\sigma,
{H}_\sigma, {\sigma})$ be the Schlichting completion. Then
$\lambda_\sigma$ is a homeomorphism from $G$ with the Hecke topology from $(G, K)$
into ${\cc B}(l^2(G, H, \sigma))$ with its weak topology.
\end{lem}

\begin{proof}
For $x, y, w\in G$ and $\xi, \eta\in V_\sigma$ we have
\begin{equation}
\bigl(\lambda_\sigma(w)\delta_{x,\xi}\mid \delta_{y,\eta} \bigr)
=\Chi_H(y^{-1}wx)\bigl(\sigma(y^{-1}wx)\xi\mid
\eta\bigr).\label{lambda_on_basis}
\end{equation}

Since $\sigma(H)$ is finite, there is $0<\varepsilon_0 <1$ such that
$K=\{h\in H\mid\, \|\sigma(h)-I\|\leq \varepsilon_0\}$. Let
$F\subseteq G$ be finite. A set of the form
$$
\mathcal{V}=\{T\in {\cc B}(l^2(G, H, \sigma))\mid \vert \bigl(
T\delta_{x,\xi_i}\mid \delta_{y,\xi_j} \bigr) -\bigl(
\delta_{x,\xi_i}\mid \delta_{y,\xi_j}\bigr)\vert<\varepsilon,
\forall x,y\in F, i,j=1,\cdots ,n\}
$$
is a typical neighbourhood of $\lambda_\sigma(e)$ for the weak
topology. We claim that $\lambda_\sigma^{-1}(\mathcal{V})=
\bigcap_{x\in F}xKx^{-1}$ if $\varepsilon\leq \varepsilon_0$. Let
$x,y\in F$ and $w\in xK x^{-1}$. Then $y^{-1}x\in H$ if and only if
$y^{-1}wx \in H$, and (\ref{lambda_on_basis}) shows that
$\lambda_\sigma(w)\in \mathcal{V}$, proving one inclusion. In
particular, $\lambda_\sigma$ is continuous at $e$, hence everywhere.
Suppose now that $\lambda_\sigma(w)\in \mathcal{V}$. Inserting $y=x$
in \eqref{lambda_on_basis} forces  $\vert
\Chi_H(x^{-1}wx)\bigl(\sigma(x^{-1}wx)\xi_i\mid
\xi_j\bigr) -\bigl(\xi_i\mid\xi_j\bigr)\vert <\varepsilon <\varepsilon_0$, so $w\in
xKx^{-1}$, and thus $\lambda_\sigma$ carries a neighbourhood subbase
at $e$ for the Hecke topology into a neighbourhood subbase at
$\lambda_\sigma(e)$ for the weak topology.
\end{proof}

Denote by $\overline{\lambda_\sigma}$ the continuous extension of
$\lambda_\sigma$ to ${G}_\sigma$. The next result shows that the induced
representation of $\sigma$ from $H_\sigma$ to $G_\sigma$ is, up to unitary equivalence,
just $\overline{\lambda_\sigma}$.

\begin{prop} With the assumptions of Lemma~\ref{induced_rep_extends}, let
$$
L^2({G}_\sigma, {H}_\sigma, {\sigma}):=\{f\in L^2({G}_\sigma,
V_\sigma)\mid f(wk)=\sigma(k^{-1}) f(w), \forall w\in {G}_\sigma,
k\in{H}_\sigma\}.
$$
Then $f\mapsto f\vert_G$ defines a unitary $U:L^2({G}_\sigma, {H}_\sigma,
{\sigma})\to l^2(G, H, \sigma)$, and $U^*\overline{\lambda_\sigma}(w) U$,
$w\in G_\sigma$, is the induced representation of $\sigma$
from $H_\sigma$ to $G_\sigma$.
\label{ind_reps_are_equiv}
\end{prop}

\begin{proof} With $(G_\sigma, K_\sigma)$ denoting the Schlichting completion of
$(G, K)$, take $w\in {G}_\sigma$ and note that since $wK_\sigma$ is
open, there is $y$ in $G$ such that $y\in wK_\sigma$. If $z\in G$,
then $w^{-1}z\in {H}_\sigma$ if and only if $y^{-1}z\in {H}_\sigma
\cap G=H$, and so $\delta_{w,\xi}=\delta_{y,\xi}$. Since
${G}_\sigma/{H}_\sigma\cong G/H$, $U$ carries the orthonormal basis
$\{\delta_{w,\xi_i} \mid w\in  {G}_\sigma /{H}_\sigma, i=1,\cdots,
n\}$ onto the orthonormal basis $\{\delta_{y,\xi_i}\mid y\in G/H,
i=1,\cdots, n\}$. Finally, a routine calculation shows that
$U^*\overline{\lambda_\sigma}(w)U$ acts as the induced
representation, and the proposition follows.
\end{proof}

With the same assumptions, let $L$ and $R$ denote the left,
respectively the right regular representation of $G_\sigma$. A
consequence of Proposition~\ref{ind_reps_are_equiv} is that
$L^2({G}_\sigma, {H}_\sigma, {\sigma})$ is a closed subspace of
$L^2(G_\sigma, V_\sigma)=L^2(G_\sigma)\otimes V_\sigma$
which is invariant under $L\otimes
I$, and that $\overline{\lambda_\sigma}$ is the restriction of
$L\otimes I$ to this subspace.

\begin{lem} Let $(G, H, \sigma)$ be a reduced Hecke triple. Then
$\widetilde R$ defined by
\begin{equation}
(\widetilde R(f)\xi)(y)=\sum_{z\in G/H}
\Delta_K(z)^{1/2} f(z)\xi(yz), \label{right_reg_rep_of_Hecke_alg}
\end{equation}
for $f\in \H_\sigma(G, H)$ and $\xi \in l^2(G, H, \sigma)$,
is a nondegenerate
$*$-representation of $\H_\sigma(G, H)$.
\end{lem}
\begin{proof}
It is straightforward to show that $\widetilde R(f)$ is well-defined and that
we get a nondegenerate
$*$-representation of $\H_\sigma(G, H)$.
\end{proof}

In \cite{Cu1}, the Hecke von Neumann algebra of $(G, H,\sigma)$ is
defined as the von Neumann algebra generated by the image of the
generalised Hecke algebra in  a left regular representation on the
space of the induced representation. Similar to this we let
$\mathcal{R}(G, H, \sigma)$ be the von Neumann algebra generated by
$\widetilde R(\H_\sigma(G, H))$ in ${\cc B}(l^2(G, H,\sigma))$. Since
Proposition~\ref{ind_reps_are_equiv} and
Proposition~\ref{Hecke_algebra_for_G_and_Gbar_same} imply that
$U\widetilde R(f)=\widetilde R(\Psi(f))U$ for all $f\in \H_\sigma({G}_\sigma,{H}_\sigma)$,
we recover an analogous result to \cite[Theorem 3]{Cu0}.

\begin{prop}
With $U$ defined in Proposition~\ref{ind_reps_are_equiv}, the map $a\mapsto
 UaU^*$
implements an isomorphism of $\mathcal{R}({G}_\sigma, {H}_\sigma,
{\sigma})$ onto $\mathcal{R}(G, H, \sigma).$
\label{Hecke_vN_same}
\end{prop}

The spaces $\mathcal{R}({G}_\sigma):=\{R_x\mid x\in {G}_\sigma\}''$
and $\mathcal{L}({G}_\sigma):=\{L_x\mid x\in {G}_\sigma\}''$ are
known to be each others commutant inside ${\cc B}(L^2({G}_\sigma))$.
Let $P_\sigma:= \int_{H_\sigma} R_k\otimes \sigma(k)\,dk$
 be the
projection corresponding to the subspace $L^2(G_\sigma,
H_\sigma,\sigma)$ of $L^2(G_\sigma, V_\sigma)$. We can now formulate
the main result of this section.

\begin{thm}
Suppose that $(G, H, \sigma)$ is a reduced Hecke triple. Then
$\lambda_\sigma(G)'$ equals $\mathcal{R}(G, H, \sigma)$.
\label{commutant_is_Hecke_vNalg}
\end{thm}

\begin{proof}
Let $({G}_\sigma, {H}_\sigma, {\sigma})$ be the
Schlichting completion of $(G, H, \sigma)$. Proposition~\ref{ind_reps_are_equiv}
implies that $\ad U^*$ carries the subset $\lambda_\sigma(G)^{''}$ of
${\cc B}(l^2(G, H,\sigma))$ into $\overline{\lambda_\sigma}(G_\sigma)^{''}$ inside
${\cc B}(L^2(G_\sigma))$. Then, since $P_\sigma$ commutes with
$L_x\otimes I$ we have
$$
\overline{\lambda_\sigma}(G_\sigma)^{''} =\{(L_x\otimes
I)P_\sigma\mid x\in G_\sigma\}^{''} =(P_\sigma \{L_x\otimes I\mid
x\in G_\sigma\}' P_\sigma)' .
$$
So by the Double Commutant Theorem we have $
\overline{\lambda_\sigma}(G_\sigma)'
=P_\sigma(\mathcal{L}(G_\sigma)'\otimes
{\cc B}(V_\sigma))P_\sigma=P_\sigma(\mathcal{R}({G}_\sigma)\otimes
{\cc B}(V_\sigma) )P_\sigma=\mathcal{R}(G_\sigma, H_\sigma, \sigma)$, and by applying
$\ad U$ we have back in ${\cc B}(l^2(G, H,\sigma))$ that $\mathcal{R}(G, H, \sigma)$
is equal to $\lambda_\sigma(G)'.$
\end{proof}

As an immediate consequence of Theorem~\ref{commutant_is_Hecke_vNalg} and
Lemma~\ref{basis_for_Hecke_algebra} we obtain the following
classical result, see \cite[Theorem 6']{M}.

\begin{cor} If $\dim\sigma=1$ then
$\lambda_\sigma$ is irreducible if and only if $B=H$.
\end{cor}

\section{$C^*$-completions of generalised Hecke algebras}\label{completions}

A consequence of Theorem~\ref{fundam_proj} is that
$p_\sigma C^*(G_\sigma,{\cc B}(V_\sigma))p_\sigma$ is
a $C^*$-completion of $\H_\sigma(G, H)$.
In this section we establish that this completion is Morita-Rieffel
equivalent to an ideal in $C^*(G_\sigma)$, and for this ideal we
describe the nondegenerate representations. We denote by
$C^*(\H_\sigma(G, H))$ the enveloping $C^*$-algebra of $\H_\sigma(G,
H)$, when it exists. We shall later be concerned with the problem of
deciding when $C^*(\H_\sigma(G, H))$ is $p_\sigma
C^*(G_\sigma)p_\sigma$.

 Note that if
$\dim \sigma=1$ we have that $p_\sigma C^*(G_\sigma)p_\sigma$ is
Morita-Rieffel equivalent to the ideal
\begin{equation}
I:=\overline{C^*(G_\sigma)p_\sigma C^*(G_\sigma)}
\label{def_of_I}
\end{equation}
in $C^*(G_\sigma)$.

We first establish that the similar $C^*$-completion of $\H_\sigma(G, H)$ can be obtained the same way if we just assume that $\sigma$ is
 finite-dimensional. Let $M$ be a compact open subgroup of
a locally compact group $L$ and $\sigma$ an irreducible
representation of $M$ on $V_\sigma$ (so in particular $V_\sigma$ is
finite-dimensional).
Then $\H_\sigma(L, M)$
is a
$*$-subalgebra of $A=C^*(L)\otimes {\cc B}(V_\sigma)$
with the convolution from \eqref{def_of_conv} and the
involution inherited from $C_c(L, {\cc B}(V_\sigma))$.
 With
$\lambda(k)$ denoting the unitary in $M(C^*(L))$ corresponding to
$k\in M$, we have that
    \[
    p_\sigma= \int_M \lambda(k)\otimes \sigma(k)\,dk
\]
in $A$, and the closure of $\H_\sigma(L, M)$ in $A$ is
$p_\sigma A p_\sigma $. If we let $d(\sigma)$ be the dimension of $V_\sigma$ and
$\operatorname{Tr}$ the unnormalised trace on
${\cc B}(V_\sigma)$, we obtain an idempotent
\begin{equation}
    \psi_\sigma(k)= d(\sigma)\operatorname{Tr}( \sigma(k))
\end{equation}
in $C^*(M)$; note that $\psi_\sigma$ is a projection in $M(C^*(L))$.
The next theorem expands upon the claims from \cite[Example
6.8]{Rie} and \cite[p 93-94]{Fe}  that Morita equivalence of
$*$-algebras gives an alternate approach to  Godement's theory of
generalised spherical functions. We do not claim  any originality in
the next result, but we rather spell out the details  in the
language of generalised Hecke algebras.

\begin{thm}\label{dim_sigma_not_one}
The generalised Hecke $C^*$-algebra $p_\sigma A p_\sigma $ is Morita-Rieffel equivalent
to the ideal $\overline{C^*(L)\psi_\sigma C^*(L)}$ in $C^*(L)$.
\end{thm}

\begin{proof}
 It is standard that
$Y=C^*(L)\otimes V_\sigma$ is a left-$A$ and right-$C^*(L)-$module in an obvious way. If we let
$T_{\xi,\eta}$ denote the rank-one operator $\alpha \mapsto (\alpha\mid\eta)\xi$,  we have a
left $A$-valued inner product on $Y$ given by
    \[
    {{_{A}}{\langle}} a\otimes \xi, b\otimes\eta \rangle= ab^*\otimes T_{\xi,\eta},
\]
and a right $C^*(L)$-valued inner product given by
    \[
    \langle a\otimes \xi, b\otimes\eta \rangle_{C^*(L)}= a^*b (\eta\mid \xi);
\]
since $V_\sigma$ is finite-dimensional, $Y$ is complete for these inner-products and so
becomes an $A-C^*(L)$ imprimitivity bimodule in Rieffel's sense. Restrict now $Y$ to the bimodule
$Y_\sigma:=p_\sigma Y$.
Clearly
\begin{align*}
_{A}\langle Y_\sigma, Y_\sigma\rangle
&=\clspn \{
    _{A}\langle p_\sigma(a\otimes \xi), p_\sigma(b\otimes\eta) \rangle\}\\
    &= p_\sigma\clspn \{_{A}\langle a\otimes \xi, b\otimes\eta \rangle \}p_\sigma\\
    &=p_\sigma A p_\sigma.
\end{align*}

On the other hand, since
 $ \clspn \{ \int \lambda( k)   (\sigma(k) \eta\mid \xi) \,
dk \}= \psi_\sigma C^*(M)$, we have
\begin{align*}
{{\langle}} Y_\sigma, Y_\sigma\rangle_{C^*(L)}
&=\clspn \{{\langle}p_\sigma(a\otimes \xi), p_\sigma(b\otimes\eta) \rangle_{C^*(L)}\}\\
&=\clspn \{\int\int
     {{\langle}}\lambda(h)a\otimes \sigma(h) \xi,
    \lambda(k)b\otimes \sigma(k) \eta \rangle_{C^*(L)} \, dh\,dk\}\\
    &= \clspn \{\int\int a^*\lambda(h\inv k)b
    (\sigma(k) \eta\mid \sigma(h) \xi) \, dh\,dk\}\\
    &= \clspn \{\int a^*\lambda( k)b
    (\sigma(k) \eta\mid \xi) \, dk\}\\
&= \clspn \{ a^* \psi_\sigma b \}=\overline{C^*(L)\psi_\sigma
C^*(L)}.
\end{align*}

Hence $Y_\sigma$ implements the equivalence between $p_\sigma A p_\sigma$ and
$\overline{C^*(L)\psi_\sigma C^*(L)}$.
\end{proof}

Thus, if $(G, H)$ is a Hecke pair and $\sigma$ is a finite-dimensional unitary representation of
$H$ with $\sigma(H)$ finite,  we form the Schlichting completion
$(G_\sigma, H_\sigma, \sigma)$ as in section~\ref{start}, and it follows from the preceeding
considerations that
$p_\sigma (C^*(G_\sigma)\otimes {\cc B}(V_\sigma))p_\sigma$
is a C$^*$-completion of $\H_\sigma(G, H)$ which is Morita-Rieffel equivalent to an ideal in
$C^*(G_\sigma)$.

\vskip 0.2cm
We now go back to the case of $1$-dimensional Hecke triples $(G, H,\sigma)$, and
we introduce the analogue of the $H$-smooth representations of $G$ that arise from $(G, H)$.

\begin{defn}
Suppose that $(G, H,\sigma)$ is a $1$-dimensional
Hecke triple.
Given a unitary representation $\pi$ of $G$ on a Hilbert space $V$,  let
\begin{equation}
V_{H, \pi}:=\{\xi \in V \mid \pi(h)\xi=\sigma(h)\xi, \forall h\in H\}.
\end{equation}
We say that $\pi$ is \emph{unitary $(H, \sigma)$-smooth} if $\clspn(\pi(G)V_{H, \pi})=V.$
\end{defn}

Note that $\lambda_\sigma$ is unitary $(H, \sigma)$-smooth (and ``smooth'' in the sense of \cite[\S 1.7]{S}), as is every representation of $G$ that is unitarily
equivalent to $\lambda_\sigma$. We have the following generalisation of
\cite[Proposition 5.18]{KLQ}.

\begin{prop}\label{smooth_rep_extends}
Let $(G, H, \sigma)$ be a reduced $1$-dimensional Hecke triple, and
$({G}_\sigma, {H}_\sigma, {\sigma})$ its Schlichting completion.
Then a representation of $G$ is unitary $(H, \sigma)$-smooth if and
only if it extends to a continuous unitary $({H}_\sigma,
{\sigma})$-smooth representation of ${G}_\sigma.$
\end{prop}

\begin{proof} The proof of \cite[Proposition 5.18]{KLQ} carries
over to this situation with one
modification: we need to show that a unitary $(H, \sigma)$-smooth representation
$\pi$ of $G$ is continuous from $G$ with the Hecke topology from
$(G, K)$ into ${\cc B}(V)$ with the strong topology. So suppose that
$x\to e$ in $G$, and pick
$\xi \in V$. By the assumption on smoothness we may take $\xi$ of the form
$\pi(y)\eta$ with $\eta\in V_{H, \pi}$. Then eventually $x$ belongs to the
neighbourhood $yK y^{-1}$, and
$$
\pi(x)\pi(y)\eta=\pi(y)\pi(y^{-1}xy)\eta=\pi(y)\sigma(y^{-1}xy)\eta=\pi(y)\eta,
$$
so $\Vert \pi(x)\pi(y)\eta - \pi(y)\eta\Vert\to 0$, proving the desired
continuity.
\end{proof}

However, the generalisation to Hecke triples of \cite[Corollary 5.10]{KLQ}
fails for $\sigma\neq 1$.

\begin{prop}
If $\sigma\neq 1$ the trivial representation of $G_\sigma$ is not unitary
$(H_\sigma, \sigma)$-smooth, and
$p_\sigma$ is not full in $C^*(G_\sigma)$.
\end{prop}

\begin{proof}
It suffices to note, first, that the integrated form of the trivial
representation  of ${G}_\sigma$ carries $p_\sigma$ into $0$, and second, that
$V_{H_\sigma, \pi}=\pi(p_{\overline{\sigma}})V$ for any representation $\pi$ of $G_\sigma$
on $V$.
\end{proof}

We show next how the strategy developed in \cite[\S 5]{KLQ}, based on Fell's
imprimitivity bimodules for $*$-algebras, for studying the
representations of $\H(G, H)$ can be carried over to tie up the
representations of $\H_\sigma(G, H)$ and the unitary
$(H, \sigma)$-smooth representations of $G$. We recall that
if ${}_{E}X_{D}$ is an imprimitivity bimodule of $*$-algebras $E$ and $D$ then a representation
$\pi$ of $D$ is \emph{positive} with respect to the right inner product
$\langle \rangle_R$ provided that $\pi(\langle f, f\rangle_R)\geq 0$ for all $f\in X$.

Let $L$ be a locally compact group, $M$ a compact open subgroup, and
$\rho$ a non-trivial character on $M$. Choose a Haar measure on $L$
such that $\rho$ becomes a projection $p_\rho$ in $C_c(L)$. Consider
the $*$-algebras $E=C_c(L)p_\rho C_c(L)$, $D=\H_\rho(L, M)$,
$B=\overline{L^1(L)p_\rho L^1(L)}$ (closure taken in $L^1(L)$) and
$C=p_\rho L^1(L)p_\rho$. (Do not confuse this $B$ with the $B$ in
\eqref{special_x}.) Then we have an inclusion of bimodules
${}_{E}X_{D}\subset {}_{B}Y_{C}$, where $X:=C_c(L)p_\rho$ and
$Y:=L^1(L)p_\rho$ with bimodule operations inherited from $L^1(L)$,
and right-inner product $\langle f, g\rangle_R=f^*g$ for $f, g$ in
$Y$ and $X$ respectively. We claim that the $C^*$-completions
$C^*(E)$ and $C^*(B)$ coincide with $\overline{C^*(L)p_\rho
C^*(L)}$. Indeed, the proof of \cite[Theorem 5.7]{KLQ} carries
through with the following alterations: given a nondegenerate
representation $\pi$ of $E$ on a Hilbert space $V$, the formula
$$
\tilde{\pi}(x)\pi(f)\xi:=\pi(xf)\xi \text{ for }x\in L, f\in E, \xi \in V,
$$
defines a representation of
$L$ on $V$ which is unitary $(M, \rho)$-smooth because $mp_{\overline{\rho}}=
\rho(m)p_{\overline{\rho}}$ for $m\in M$ implies that
$\pi(p_{\overline{\rho}})V= V_{M, \rho}$. The integrated form
of $\tilde{\pi}$ will then be a nondegenerate extension of $\pi$ to
$\overline{C^*(L)p_{\overline{\rho}} C^*(L)}$, from which the claim follows. So we
obtain the analogue of the category equivalences from \cite[Corollaries 5.12 and 5.20]{KLQ}.

\begin{cor} Let $(G, H,\sigma)$ be a reduced $1$-dimensional Hecke triple. Then there are
category equivalences between

\textnormal{a)} the unitary $(H,\sigma)$-smooth representations of $G$ and the
$\langle \rangle_R$-positive representations of $\H_\sigma(G, H)$,

\textnormal{b)} the nondegenerate representations of
$\overline{C^*(G_\sigma)p_\sigma
C^*(G_\sigma)}$ and the $\langle \rangle_R$-positive representations of
$p_\sigma L^1(G_\sigma)p_\sigma.$
\label{category_equiv}
\end{cor}

\section{The case $H\subseteq N\trianglelefteq G$.}\label{H_in_N_normal}

Let $(G, H, \sigma)$ be a reduced $1$-dimensional Hecke triple with
Schlichting completion $(G_\sigma, H_\sigma, \sigma)$, and suppose
that $H$ is contained in a normal subgroup $N$ of $G$. (This is
clearly interesting only if $N\neq G$.) We will show that $I$
defined in \eqref{def_of_I} is a (twisted) crossed product, see
\cite{GKP} (and \cite{PaR}) for definitions. Let $N_\sigma$ denote
the closure of $N$ in $G_\sigma$, and let $\ad$ be the action of
$G_\sigma$ by conjugation on $N_\sigma$ (and on $C^*(N_\sigma)$).
The universal covariant representation $(\pi, u)$ of
$(C^*(N_\sigma), G_\sigma)$ into the twisted crossed product
$C^*(N_\sigma)\rtimes G_\sigma/N_\sigma$ determines an isomorphism
$\pi\times u:\pi(b)u(f)\mapsto bf$ of $C^*(N_\sigma)\rtimes
G_\sigma/N_\sigma$ onto $C^*(G_\sigma)$. Note that the twist
disappears when $G$ is a  semi-direct product of $N$ by a group $Q$,
see also \cite{LarR}.

\begin{thm} Suppose that $(G, H, \sigma)$ is a reduced $1$-dimensional Hecke triple such that
$H$ is contained in a normal subgroup $N$ of $G$. Let $N_\sigma$
denote the closure of $N$ in the Schlichting completion $(G_\sigma,
H_\sigma, \sigma)$ of $(G, H, \sigma)$.

\smallskip
\textnormal{(a)} Then
$
I_\sigma:=\clspn\{xp_\sigma x^{-1}n\mid x\in G_\sigma, n\in N_\sigma\}
$
is an $\ad$-invariant ideal of $C^*({N}_\sigma)$, and the isomorphism
$\pi\times u$ carries $I_\sigma\rtimes G_\sigma/N_\sigma$  onto $I$ defined in
\eqref{def_of_I}.

\smallskip
\textnormal{(b)} $I_\sigma \rtimes G/N$ is Morita-Rieffel equivalent
to $p_\sigma C^*(G_\sigma)p_\sigma$.
\label{Hecke_alg_as_cp}
\end{thm}

\begin{proof}
Since $xH\mapsto xH_\sigma$ is a bijection from $G/H$ onto
$G_\sigma/H_\sigma$ (essentially by \cite[Proposition 3.9]{KLQ}), it follows that
$$
\spn\{xp_\sigma x^{-1}n\mid x\in G_\sigma, n\in N_\sigma\}=
\spn\{xp_\sigma x^{-1}n\mid x\in G, n\in N\}.
$$
Note that $xp_\sigma x^{-1}\in C_c(N_\sigma)$, so $I_\sigma \subset C^*(N_\sigma)$.  Since $mI_\sigma=I_\sigma m=
I_\sigma$ for $m\in N$, we get as in the proof of \cite[Theorem 8.1]{KLQ}
that $I_\sigma$ is a closed, $\ad$-invariant ideal
of $C^*(N_\sigma)$. Since $xp_\sigma y=xp_\sigma x^{-1}(xp_\sigma y)\subset
\spn\{xp_\sigma x^{-1}f\mid f\in C_c(G_\sigma)\}$, it follows that
\begin{align*}
(\pi\times u)(I_\sigma \rtimes G_\sigma/N_\sigma)
&=\clspn\{xp_\sigma x^{-1}nf\mid x\in G, n\in N, f\in C_c(G_\sigma)\}\\
&=\clspn\{xp_\sigma x^{-1}f\mid x\in G, f\in C_c(G_\sigma)\}\\
&=\clspn\{xp_\sigma y\mid x, y\in G\}=I,
\end{align*}
as claimed in (a).

For (b), it suffices by (a) to establish that $I_\sigma \rtimes G/N\cong
I_\sigma \rtimes G_\sigma/N_\sigma$. Since $hp_\sigma=\overline{\sigma(h)}
p_\sigma$ for all $h\in H$, the argument in the proof of \cite[Theorem 8.2]{KLQ}
shows that the canonical
homomorphism $\omega:G\to M(I_\sigma \rtimes G/N)$ is  unitary $(H, \sigma)$-smooth.
Thus  $\omega$ has a continuous extension
$\overline{\omega}$ to $G_\sigma$ by Proposition~\ref{smooth_rep_extends}, and
then $\overline{\omega}$ forms a covariant
pair together with the canonical homomorphism $I_\sigma \to M(I_\sigma \rtimes
G/N)$, from which the claimed isomorphism follows.
\end{proof}

\begin{cor} With the notation from Theorem~\ref{Hecke_alg_as_cp},
assume that $B$ is a subgroup of $G$ such that $N$ is a normal subgroup
of $B$. Then $I_{\sigma, B}=\clspn\{xp_\sigma x^{-1}n\mid x\in B_\sigma, n\in
N_\sigma\}$ is an $\ad$-invariant ideal of $C^*(N_\sigma)$, and the closed ideal
generated by $p_\sigma$ in $C^*(B_\sigma)$ is Morita-Rieffel equivalent
to the twisted crossed products $I_{\sigma, B}\rtimes B_\sigma/N_\sigma$
and $I_{\sigma, B}\rtimes B/N$.
\label{B_Hecke_alg_as_cp}
\end{cor}

With the hypotheses of Theorem~\ref{Hecke_alg_as_cp}, we now assume that $N$ is abelian, and we consider  the Fourier transform $f\mapsto \hat{f}$
from $C^*({N}_\sigma)$ onto $C_0(\widehat{N_\sigma})$. We let
\begin{equation}
\sigma +H_\sigma^\perp:=\{\alpha\in \widehat{N_\sigma}\mid
\alpha\vert_{H_\sigma}={\sigma}\}
\label{def_of_H_perp}
\end{equation}
be the set of all continuous extensions of $\sigma$ to $N_\sigma$. One can
verify that $\widehat{p_{\sigma}}=\Chi_{\sigma +H_\sigma^\perp}$.
The dual action of $G_\sigma$ on $\widehat{N_\sigma}$ is characterised by
$$
\langle n, x \cdot\alpha\rangle=\langle x^{-1} nx, \alpha\rangle
\text{ for }n\in N_\sigma, \alpha \in \widehat{N_\sigma}, x\in G_\sigma,
$$
and then we have  for all $x\in G$ that
\begin{align}
(xp_\sigma x^{-1})^\wedge(\alpha)
&=\int_{N_\sigma}\overline{\langle m ,\alpha\rangle}(xp_\sigma x^{-1})(m)d\mu(m)
\notag\\
&=\Delta(x)\int_{N_\sigma}\overline{\langle m, \alpha\rangle}\langle
x^{-1}mx, \sigma \rangle\Chi_{xH_\sigma
x^{-1}}(m)d\mu(m)\notag\\
&=\Delta(x)\int_{xH_\sigma x^{-1}} \langle m, x\cdot \sigma-\alpha\rangle d\mu(m)
\notag\\
&= \Delta(x)\mu(xH_\sigma x^{-1})\Chi_{x\cdot(\sigma +H_\sigma^\perp)}(\alpha)\notag
\\
&=\Chi_{x\cdot(\sigma +H_\sigma^\perp)}(\alpha).
\label{Ft_of_p_sigmax}
\end{align}
Therefore $\widehat{I_\sigma}$ is the ideal in $C_0(\widehat{N_\sigma})$
generated by $\{ \Chi_{x\cdot(\sigma +H_\sigma^\perp)}\mid x\in G\}$, so if we
let
\begin{equation}
\Omega_{\sigma}:=\bigcup_{x\in G} {x}\cdot(\sigma +H_\sigma^\perp),
\label{def_omega}
\end{equation}
then  we have proved the following specialisation of Theorem~\ref{Hecke_alg_as_cp}:

\begin{thm}
With the hypotheses of Theorem~\ref{Hecke_alg_as_cp}, if $N$ is moreover
abelian, then  $\widehat{I_\sigma}=C_0(\Omega_{\sigma})$. If $\sigma$ is
non-trivial, then
$0\notin \Omega_{\sigma}$, so $p_\sigma$ is not full.
\label{theo_N_abelian}
\end{thm}

\begin{rem} If $B$ is a group and $N$ is normal in $B$, then
$\widehat{I_{\sigma, B}}$ from Corollary~\ref{B_Hecke_alg_as_cp}
equals $C_0(\Omega_{\sigma, B})$, where
\begin{equation}
\Omega_{\sigma, B}:=\bigcup_{x\in B} x\cdot (\sigma + H_\sigma^\perp).
\label{def_omegaB}
\end{equation}
\end{rem}

\subsection{The case $H\trianglelefteq N\trianglelefteq G$.}
\label{H_normal_N_normal}
Suppose that $(G, H)$ is a reduced Hecke pair and $N$ is such that $H\trianglelefteq N\trianglelefteq G$. Let $\sigma$ be a finite character of $H$ such that
\begin{equation}
\sigma(nhn^{-1})=\sigma(h) \text{ for
all }n\in N, h\in H,
\label{N_in_B}
\end{equation}
(so $N\subseteq B$
in the notation of \eqref{special_x}). Then $p_\sigma$ is central in $C^*(N_\sigma)$,
and so are $xp_\sigma x^{-1}$ for all $x$ in $G_\sigma$. On one hand, this shows that
the ideal $I_\sigma$ defined in
Theorem~\ref{Hecke_alg_as_cp} will satisfy $I_\sigma=C^*(N_\sigma)p_1$, where
$p_1:=\operatorname{sup}\{xp_\sigma x^{-1}\mid x\in G_\sigma\}$.

On the other hand, $xp_\sigma x^{-1}p_\sigma$ is a projection, and
$p_\sigma xp_\sigma$ is a partial
isometry, for all $x\in G_\sigma$. Thus for every $*$-representation $\pi$ of
$\H_\sigma(G, H)$ and for all $x\in B$ we have
$$
\Vert \pi(p_\sigma x p_\sigma)\Vert\leq 1=\Vert p_\sigma x p_\sigma\Vert_1,
$$
where the equality is from Lemma~\ref{psigma_x_psigma}. This establishes that
 $C^*(\H_\sigma(G, H))$ exists and is equal to the enveloping
$C^*$-algebra of $l^1(G, H,\sigma)$.

But more is true. We show next that the
right inner product $\langle f,g \rangle_R=f^*g$ on $X:=C_c(G_\sigma)p_\sigma$ is
positive  in the following sense: given $f$ in $X$,
there are $g_i$ in $\H_\sigma(G_\sigma, H_\sigma)$, $i=1,\dots ,n$, such that
$\langle f, f\rangle_R=\sum_{i=1}^n g_i^*g_i$.

\begin{thm}
Let $L$ be a locally compact totally disconnected group, $M$ a compact open
subgroup, and $\rho$ a
non-trivial character on $M$. Suppose
that $M$ is normal in a closed normal subgroup $N$ of $L$, and choose a Haar
measure on $L$ such that $p(m)=\rho(m)\Chi_M(m)$ becomes a self-adjoint projection in $C_c(L)$.
If
$$
\rho(nmn^{-1})=\rho(m) \text{ for all }n\in N, m\in M,
$$
then $Y:=C_c(L)p$ is a left-$C_c(L)p C_c(L)$ and right-$\H_\rho(L, M)$ bimodule
of $*$-algebras with positive right inner product. Moreover,
$C^*(\H_\rho(L, M))=p C^*(L)p.$
\label{C*_completion_normal_subgroups_general}
\end{thm}

We have the following consequence of this theorem and of
Proposition~\ref{H_sigma_spanned_by_epsilonx}.

\begin{cor}\label{C*_completion_normal_subgroups}
 Let $(G, H, \sigma)$ be a reduced $1$-dimensional  Hecke triple, $N$ a normal subgroup of $G$
such that $H$ is normal in $N$, and suppose that $\sigma$ satisfies
\eqref{N_in_B}. Let $(G_\sigma, H_\sigma, \sigma)$ be the
Schlichting completion of $(G, H, \sigma)$. Then
$C^*(\H_\sigma(G_\sigma, H_\sigma))=p_\sigma C^*(G_\sigma)p_\sigma.$

If $B$ is a group, then $C^*(\H_\sigma(G_\sigma,
H_\sigma))=p_\sigma C^*(B_\sigma) p_\sigma$.
\end{cor}

\begin{proof}[Proof of Theorem~\ref{C*_completion_normal_subgroups_general}]
The bimodule $Y$ is spanned by $\{xp\mid x\in L\}$. By the hypothesis on $\rho$,  the
projection $p$ is central in $C_c(N)$ and hence, by normality of $N$ in
$L$, so is $xp x^{-1}$ for every $x\in L$. Then
$\{xp x^{-1}\mid x\in L\}$ are commuting projections in $C^*(L)$ and by
following verbatim the proof of  \cite[Theorem 5.13]{KLQ} we
conclude that for any element $f=\sum_1^n c_i x_ip$ in $Y$, the product
$f^*f$ is a finite sum of elements $h^*h$ with $h\in p C_c(L)p$. The last claim
follows from \cite[Proposition 5.5(iii)]{KLQ} applied to the bimodule $Y$.
\end{proof}

\section{The case when $(B, H)$ is directed.}

Let $(G, H,\sigma)$ be a  reduced $1$-dimensional
 Hecke triple with Schlichting completion $(G_\sigma, H_\sigma, \sigma)$.
Let $B$ be the subset of $G$ defined in \eqref{special_x}, and consider its
closure $B_\sigma$ from Proposition~\ref{Hx_Kx}.

\begin{lem} If $x\in B$ and $xHx^{-1}\supset H$, then
$x^{-1}p_\sigma x\geq p_\sigma$ in $C^*(G_\sigma)$.
\label{B_directed}
\end{lem}

\begin{proof}
It suffices to note that $x^{-1}p_\sigma x=\mu(x^{-1}H_\sigma x)^{-1}\sigma
\vert_{x^{-1}H_\sigma x}.$
\end{proof}

Assume  that $B$ is a group. Consider the semigroup
$B^{+}:=\{x\in B\mid xHx^{-1}\supset H\}$, and recall from
\cite[Definition 6.1]{KLQ} that $(B, H)$ is called  directed if $B^{+}$ is an Ore
semigroup in $B$, i.e. $B=\{x^{-1}y\mid x,y \in B^{+}\}$. In general, $(B, H)$ need not be
directed when $(G, H)$ is.
For $x,y$ in $B$ with $yx\in B^{+}$ we have
$x^{-1}y^{-1}p_\sigma yx\geq p_\sigma$ by Lemma~\ref{B_directed},
so $p_\sigma yxp_\sigma=yxp_\sigma$, and the proof of \cite[Theorem 6.4]{KLQ} can be used here to give that the bimodule $C_c(B_\sigma)p_\sigma$
has
positive right inner product. Therefore \cite[Proposition 5.5 (iii)]{KLQ} gives the following.

\begin{prop}
Let $(G, H, \sigma)$ be a reduced $1$-dimensional
 Hecke triple such that $B$ is a subgroup of
$G$ and $(B, H)$ is directed. Then
$
C^*(\H_\sigma(G, H))=C^*(\H_\sigma(B, H))=p_\sigma C^*(B_\sigma)p_\sigma.
$
\label{C*_completion_B_directed}
\end{prop}

In fact we have a precise description of the ideal in $C^*(B_\sigma)$
generated by $p_\sigma$.

\begin{thm}
Let $(G, H, \sigma)$ be a reduced $1$-dimensional
 Hecke triple such that $B$ is a subgroup of
$G$ and $(B, H)$ is directed. Denote by $H_\infty$ the subgroup
$\bigcap_{x\in B^+} x\inv H_{\sigma}x$ of $H_\sigma$,
and let $\mu_\infty$ be  normalised Haar measure on $H_\infty$. Then the closed
ideal generated by $p_\sigma$ in $C^*(B_\sigma)$ is equal to
$C^*(B_\sigma)p_{\sigma,\infty}$, where
\begin{equation}
p_{\sigma,\infty}=\int_{H_\infty}{\sigma(h)}hd\mu_\infty(h).
\label{def_of_p1}
\end{equation}
\label{theo_p1_B_directed}
\end{thm}

To prove this theorem we need a general lemma.

\begin{lem}
Suppose that $G$ is a locally compact group,
$H$ a compact open subgroup, $\{H_i\}_{i\in I}$  a family of open subgroups
of $H$ over a directed set $I$ with $H_j\subset H_i$ for $i<j$, and $\sigma$
a character of $H$ with finite range. Denote $H_\infty=\cap_{i\in I}H_i$,
let $\mu_i$ be normalised Haar measure
on $H_i$ for $i\in I\cup \{\infty\}$, and let
$$
p_{\sigma,i}=\int_{H_i}{\sigma(h)}hd\mu_i(h)
$$
for $i\in I\cup \{\infty\}$. Then $p_{\sigma,i} \rightarrow
p_{\sigma,\infty}$ in $M(C^*(G))$.
\label{last_small_lemma}
\end{lem}

\begin{proof}
Let $K_i=H_i\cap \ker\sigma$ and $\nu_i$ be normalised Haar
measure on $K_i$ for $i\in I\cup \{\infty\}$. Define $q_{i}=\int_{K_i}hd\nu_i(h)$ for $i\in I\cup
\{\infty\}$.

We claim first that $q_{i} \rightarrow q_{\infty}$ in $M(C^*(G))$.
Given $a\in C^*(G)$ and $\epsilon>0$, take $b=q_{\infty}a$ and let $
U=\{x\in G\mid \|xb-b\| <\epsilon\}$.
Then $U$ is open and contains $K_\infty$. We want to show that $K_i\subset
U$ eventually. If not, for every $i\in I$ there is $k_i\in K_i\setminus U$,
and then by compactness of $H$ there is a subnet $(k_i)$ converging to an
element $k\in H\setminus U$. For given $i$, the set $K_ik$ is open and contains $k$,
and thus there is $j>i$ such that $k_j\in K_ik$.
Therefore $k\in K_ik_j\subset K_iK_j=K_i$, and it follows that
$k\in \cap K_i=K_\infty$, contradicting $k\notin U$. Since $q_iq_\infty =q_i$ we have
\[
q_ia-q_\infty a =q_ib-b=\int_{K_i} (hb-b)\,d\nu_i(h).
\]
But  $K_i\subset U$ eventually, and so
$\| q_ia-q_\infty a\|<\epsilon$ for sufficiently large $i$, proving the claim.

We next claim that $p_{\sigma,i}=q_{i} p_{\sigma,\infty}$ for sufficiently large $i$.
Indeed, since the finite sets $\sigma(H_i)$ form a decreasing family,
there is $i_0\in I$ such that $\sigma(H_i)=\sigma(H_{i_0})$ for $i>i_0$. Then
$\sigma(H_\infty)=\cap_j \sigma(H_{j})=\sigma(H_{i_0})=\sigma(H_{i})$ for $i>i_0$,
and hence $H_i=K_iH_\infty$ $i>i_0$. Then for each $i>i_0$,
the Haar measure $\mu_i$ on $H_i$ is given by
\[
\int_{H_i}f(h)\, d\mu_i(h)
=\int_{K_i}\int_{H_\infty} f(kl) \, d\nu_i(k)\, d\mu_\infty(l),
\]
and the claim follows because
$$
p_{\sigma,i}
=\int_{H_\infty}\int_{K_i} k{\sigma(l)}l\, d\nu_i(k)\, d\mu_\infty(l)
 =q_{i} p_{\sigma,\infty}.
 $$
Using the two claims we conclude the proof of the lemma by observing that
$$
p_{\sigma,i}=q_{i} p_{\sigma,\infty}\rightarrow q_\infty p_{\sigma,\infty}=p_{\sigma,\infty}\text{ in }
M(C^*(G)).
$$
\end{proof}

\begin{proof}[Proof of Theorem~\ref{theo_p1_B_directed}.]
Since $(B,H)$ is directed, Lemma~\ref{B_directed} shows that
$x^{-1}p_\sigma x\geq p_\sigma$ for $x\in B^{+}$.
The projections $p_{\sigma, x}:=x\inv p_\sigma x$ have support in $x^{-1}H_\sigma
x$, and then Lemma~\ref{last_small_lemma} implies that $p_{\sigma, x}\nearrow
p_{\sigma, \infty}$ in $M(C^*(B_\sigma))$. Hence the ideal generated by $p_\sigma$
in
$C^*(B_\sigma)$ is $C^*(B_\sigma)p_{\sigma, \infty}$.
\end{proof}

Suppose in addition that $H$ is contained in an abelian subgroup $N$ which is normal in $B$. Then $x^{-1}p_\sigma x$ is a projection in
$C^*(N_\sigma)$ for all $x\in B^+$, and  the Fourier transform applied
to both sides of the inequality $x^{-1}p_\sigma x\geq p_\sigma$ gives
$\widehat{p_{\sigma, x^{-1}}}\geq \widehat{p_\sigma}$. Hence \eqref{Ft_of_p_sigmax}
implies that $x\cdot(\sigma +H_\sigma^\perp)\subset \sigma +H_\sigma^\perp$ for
$x\in B^+$. Since $\widehat{p_{\sigma, x^{-1}}}$ converges in
$C_0(\widehat{N_\sigma})$ to $\widehat{p_{\sigma,\infty}}$, and since
$\widehat{p_{\sigma,\infty}}$ equals the characteristic function of the set
$$
\sigma +H_\infty^\perp:=\{\alpha\in \widehat{N_\sigma}\mid \alpha\vert_{H_\infty}=
\sigma\vert_{H_\infty}\},
$$
we obtain the following strengthening of Corollary~\ref{B_Hecke_alg_as_cp} and
\eqref{def_omegaB}:

\begin{cor}
With the notation above, the Hecke algebra
$C^*(\H_\sigma(G,H))$
is Morita-Rieffel equivalent to
$$
\overline{C^*(B_\sigma) p_\sigma C^*(B_\sigma)}\cong C^*(N_\sigma)p_{\sigma, \infty}
\rtimes B/N\cong C_0(\sigma +H_\infty^\perp)\rtimes B/N.
$$
\label{B_directed_N_abelian}
\end{cor}

\section{Applications}\label{applications}

All the examples have $\dim\sigma=1$ and we begin this section by analysing a simple situation where the generalised Hecke algebra
of a reduced Hecke triple $(G, H, \sigma)$ does not depend on $\sigma$.

\begin{prop}\label{sigma_extends}
Let $(G, H, \sigma)$ be a
$1$-dimensional
Hecke triple, and suppose that $\sigma$ extends
to a character of $G$. Then the map
$$
\Phi(f)(x)=\overline{\sigma(x)}f(x)
$$
for $f\in \H_\sigma(G, H)$ and $x\in G$, is a $*$-isomorphism of $\H_\sigma(G, H)$
onto $\H(G, H).$
\end{prop}

\begin{proof} The definition of $\Phi$ implies that $\Phi(f)(hxk)=\Phi(f)(x)$
for $f\in \H_\sigma(G, H)$, $x\in G$ and $h, k\in H$, so $\Phi(f)$ is $H$-biinvariant
and thus $\Phi$ is well-defined. Using that $\Delta_H=\Delta_K$ shows  that
$$
\Phi(f)^*(x)
=\Delta_H(x^{-1})\overline{\Phi(f)(x^{-1})}=\overline{\sigma(x)}
\Delta_K(x^{-1})\overline{f(x^{-1})}=\Phi(f^*)(x),
$$
so $\Phi$ is adjoint preserving. Finally, since $\sigma$ is defined everywhere on $G$,
a routine verification shows that $\Phi(f\ast g)=\Phi(f)\ast \Phi(g)$ for all
$f,g\in \H_\sigma(G, H)$, as wanted.
\end{proof}

Criteria for when a unitary representation of $H$ extends to a unitary representation of $G$ on the same Hilbert space
have been recently analysed in, for example, \cite{anHKR}. However, in the
examples where we have been able to extend, it has been straightforward to write
down a formula for the extended character.

\begin{ex}\label{example1}
Suppose that $H$ and $N$ are subgroups of a group $G$ such that
$H$ is finite, $N\trianglelefteq G$, $G=HN$, and $(G, H)$ is reduced. Suppose
that $\sigma$
is a character on $H$. Then $\sigma(hn):=\sigma(h)$ for $h\in H$ and
$n\in N$ is a well-defined extension to a character on $G$.  Hence
 $\H_\sigma(G, H)$ is isomorphic to $\H(G, H)$ by Proposition~\ref{sigma_extends}.

As a concrete example of this set-up we can take
 $G$ to  be the infinite dihedral group $\Z\rtimes_\psi \Z_2$ with generators
 $a$ for $\Z$ and $b$ for  $\Z_2$, where $\psi_b(a)=a^{-1}$.
Alternatively, $G$ has presentation
$\langle a,b\mid b^2=1, bab=a^{-1} \rangle$. Let $H=\langle b\rangle\cong
\Z_2$, $N=\Z$, and consider the character $\sigma:H\to \T$, $\sigma(b)=-1$.
Using either \cite[Example 3.4]{T} or \cite[Example 10.1]{KLQ}
we  conclude that $\H_\sigma(G, H)$ does not
have a largest $C^*$-norm because  $\H(G, H)$ fails to have one.
\end{ex}

Suppose that $(G, H, \sigma)$ is a reduced Hecke triple. Let $(G_\sigma,
H_\sigma, \sigma)$ and $(G_0, H_0)$ be the Schlichting completions of $(G, H,
\sigma)$ and $(G, H)$, respectively. Choose left invariant Haar measures
$\mu$ on $G_\sigma$ and
$\nu$ on $G_0$,  normalised so that $\mu(H_\sigma)=1$ and $\nu(H_0)=1$. Let
$p_\sigma$ be the projection defined in Theorem~\ref {fundam_proj}, and
$p_0$ the projection $\Chi_{H_0}$ in $C_c(G_0)$. In certain cases we can
identify $p_\sigma C^*(G_\sigma)p_\sigma$
and $p_0C^*(G_0)p_0$ inside the same algebra, as shown in the next proposition.
In concrete examples, it suffices to verify whether $K\supseteq x_0Hx_0^{-1}$
for some $x_0$ in $G$, because then the continuity hypothesis is automatic.

\begin{prop}\label{Gsigma_is_G0}
Given a reduced $1$-dimensional
 Hecke triple $(G, H,\sigma)$, suppose that $\sigma$ extends to a character
of $G$, and is continuous with respect to the Hecke topology from
$(G, H)$. Then $\iota:G_\sigma\to G_0$ is a topological isomorphism,
and $\Phi(x)=\overline{\sigma(x)}x$ for $x\in G$
extends to an automorphism of $C^*(G_\sigma)$ which carries $p_\sigma$ into
$p_0$.
\end{prop}

\begin{proof}
The hypothesis implies that $\sigma$ has a continuous extension $\sigma_0$ to
$H_0$. By the second part of Theorem~\ref{uniqueness_univ_prop}, the map
$\iota$ is a topological isomorphism of $G_\sigma$ onto $G_0$ and of $H_\sigma$
onto $H_0$. Since the modular functions of $G_0$ and $G_\sigma$ coincide on $G$
by Corollary~\ref{modular_functions_same}, the involution is preserved by $\Phi$.
\end{proof}

\begin{ex}\label{rational_Heisenberg}
 (The rational Heisenberg group.) We analyse now the Hecke pair
studied in  \cite[Example 10.7]{KLQ}. We use the same notation, so
\[
[u,v,w]:=
\begin{pmatrix}
1&v&w\\
0&1&u\\
0&0&1
\end{pmatrix}, \text{ where }u, v, w\in \Q.
\]
Then $G=\bigl\{\,[u,v,w] \bigm| u,v\in \Q, w\in \Q/\Z\,\bigr\}$ and
$H=\bigl\{\, [u,v,0] \bigm| u,v\in \Z\,\bigr\}$ form a
reduced Hecke pair. We let $N$ be the (abelian) subgroup of $G$
with $u, v\in \Z$, and then  $H\trianglelefteq N \trianglelefteq G$.
Fix $s, t$ in $\Q$ and let $\sigma$ be the character of $H$ given by
\begin{equation}
\sigma([m,n,0])=\exp(2\pi i (sm +tn))\text{ for }m, n\in \Z.
\label{char_on_Heisenberg_subgroup}
\end{equation}

Since $[u_1,v_1,w_1][u_2,v_2,w_2]=[u_1+u_2,v_1+v_2,w_1+w_2+v_1u_2]$ in $G$,
the equation \eqref{char_on_Heisenberg_subgroup} extends to a character
$\sigma$ on $G$ given by $\sigma([u,v,w])=\exp(2\pi i(su+tv))$. Thus
$\H_\sigma(G, H)$ is isomorphic to $\H(G, H)$ by Proposition~\ref{sigma_extends}.
We know from \cite[Example 10.7]{KLQ} that the collection of sets
$$
H_{x,y}=\{[u,v,0]\mid u\in \Z\cap y\Z,v\in \Z\cap x\Z\}
$$
forms a neighbourhood base at $e$ when $x,y\in \Z\setminus \{0\}$.
If we denote by $b$ and $d$ the denominators of $s$ and $t$ respectively, then
$$
K=\{[m, n,0]\mid sm+tn\in \Z\}\supset H_{d, b}.
$$
Thus $\sigma$ is continuous at $e$, hence everywhere,
for  the Hecke topology  from $(G,H)$, and so Proposition~\ref{Gsigma_is_G0} implies
that the Schlichting completion $(G_0, H_0)$ of $(G, H)$ is also
the Schlichting completion of $(G, H,\sigma)$. Then
Corollary~\ref{C*_completion_normal_subgroups} or \cite[Theorem 5.13]{KLQ} imply
that $p_0C^*(G_0)p_0$ is the largest $C^*$-completion of $\H_\sigma(G, H)$.

For different  choices of $\sigma$, the ideals $\widehat{I_\sigma}\rtimes G/N$ from
Theorem~\ref{Hecke_alg_as_cp} are all isomorphic to $\widehat{I_0}\rtimes G/N$,
 where $I_0$ corresponds to the trivial character $\sigma\equiv 1$. Let
$\mathcal{A}_f$ and $\mathcal{Z}$ respectively denote the ring of finite adeles
and its compact  open subring of integral adeles. To describe the
sets $\Omega_\sigma$ defined in \eqref{def_omega}, we recall from
\cite{KLQ} that $G_0=\{[u,v,w]\mid u,v\in \mathcal{A}_f, w\in Q/Z\}$, with
$N_0$ the subgroup with components $u,v\in \mathcal{Z}$, and hereby with
$H_0=\{[u,v,0]\mid u, v\in \mathcal{Z}\}$. Then a computation shows that
$\Omega_\sigma=[s,t,0]+\Omega_0$, where $\Omega_0$ corresponding to the trivial
character $\sigma$ is described in \cite[Example 10.7]{KLQ}. The
isomorphism $\Phi$ is obtained from translation by $[s,t,0]$.
\end{ex}

\begin{ex}\label{ax_p_adic}
(The $p$-adic $ax+b$-group.)
Let
$p$ be a prime and denote
$$
N:=\Z[ p^{-1}]
=\{\frac m{p^n}\mid m, n\in \Z, n\geq 0\}.
$$
Let $G:=N\rtimes\Z$ with $m\cdot b=bp^m$ for $m\in \Z$ and $b\in N$. Then $G$
and $H:=\Z$ form a reduced Hecke pair. Let $q$  be a positive non-zero integer that is
co-prime with $p$, and $\sigma$ the character of $H$
given by $\sigma(n)=\operatorname{exp}(\frac {2\pi i n}{q})$. Then $K=\ker \sigma=q\Z$, and $(G, K)$ is also
reduced. For $g=(x, p^k)\in G$ we have that
$gHg^{-1}=p^{-k}\Z$, and hence deduce that $\sigma$ is not continuous in the
Hecke topology from $(G, H).$ We will study $\H_\sigma(G, H)$ using Corollary~\ref{B_directed_N_abelian}.

In order to describe the Schlichting completion of $(G, H, \sigma)$ we
recall some facts about ${\mathbf q}$-adic integers and numbers, where $\mathbf{q}$ is a doubly infinite sequence of integers greater than one. We refer to
\cite[\S 12.3.35]{P} for details (see also \cite[\S 25.1]{HR}).
We are interested in the particular
sequence  ${\mathbf q}$ in which ${\mathbf q}_0=q$ and ${\mathbf q}_n=p$ for
all other $n\in \Z$, and we view an element $a$ of $\Omega_{\mathbf q}$ as
a formal sum
$$
a=\sum_{i=-i_0}^{-1}a_ip^i + a_* +q\sum_{i=0}^\infty a_ip^i
$$
with $0\leq a_i<p$ and $0\leq a_*< q$. By denoting
$a_{-}:=\sum_{i=-i_0}^{-1}a_ip^i$ and $a_+:=\sum_{i=0}^\infty a_ip^i$ we have
\begin{equation}a=a_{-}+a_*+qa_+
\label{formal_a}
\end{equation}
 with $a_+\in \Z_p$. Elements $a$ in $\Omega_{\mathbf q}^0$ are characterised by the condition that
 $a_{-}=0$. The
first theorem in
\cite[\S 12.3.35]{P} says that there is an injective group homomorphism $\phi$ from $\Z[p^{-1}]$
into the locally compact, totally disconnected (additive) abelian group $\Omega_{\mathbf q}$, such that $\phi$ has
dense range, and restricts to a bijection of $\Z$ onto a dense subgroup of the compact, totally disconnected subgroup $\Omega_{\mathbf q}^0$
of $\Omega_{\mathbf q}$.

Multiplication by $p$ (with carry-over) is a continuous action of $\Z$ on our $\mathbf q$-adic numbers (because it is clearly continuous on the dense subset of
elements $a$ with $a_+$ finite), and so $\phi$ extends to a group homomorphism from $G$ into
the locally compact, totally disconnected group
$L:=\Omega_{\mathbf q}\rtimes \Z$. The range of $\phi$ is still dense, and
$M:=\Omega_{\mathbf q}^0$ is compact, open in $L$. If $y\in G$ is such that $\phi(y)\in M$, then $\phi(y)_+$ is finite and $\phi(y)_{-}=0$, and it follows that
$\phi^{-1}(M)=H$. The formula
$$
\rho(a_* +qa_+):=\operatorname{exp}(2\pi i a_*/q)
$$
defines a character of $M$. Note that the kernel $K$ of $\rho$ consists of the formal sums
$\{qa_+\mid a_+\in \Z_p\}$. Now $(L, K)$ being reduced is the same as
$\bigcap_{m\in \Z}p^mK=\{0\}$, and this last identity can be verified directly
using  \eqref{formal_a}. From the definition we have
$\rho\circ \phi\vert_H=\sigma$, and Theorem~\ref{uniqueness_univ_prop} gives
the following.

\begin{prop}
The Schlichting completion of $(G, H, \sigma)$ is
$(\Omega_{\mathbf q}\rtimes \Z, \Omega_{\mathbf q}^0, \rho)$.
\end{prop}

It is also possible to identify $\Omega_{\mathbf q}$ as the topological limit
$\varprojlim N/qp^n \Z$, where the bonding maps are reductions modulo $qp^n\Z$, see
for example \cite[Proposition 3.10]{KLQ};
then $\Omega_{\mathbf q}^0$ is the profinite group $\varprojlim \Z/qp^n\Z$.

The Schlichting completion of $(G, H)$ is $(\Q_p\rtimes \Z, \Z_p)$,
where $\Q_p$ and $\Z_p$ are respectively the $p$-adic numbers and
$p$-adic integers. Let $(a,k)\in N_\sigma\rtimes \Z$ with $a$ as in
\eqref{formal_a}. Then the homomorphism $\iota$ from
Remark~\ref{def_of_iota} sends $(a_{-}+a_*+qa_+,k)$ into $(a_{-}+a_*+a_+, k)$,
and so is not a topological isomorphism.

\begin{lem}\label{B_for_BS}
 Let $n_0$ be the smallest integer  $n>0$ such that $p^n\equiv
1\mod q$. Then $B=\Z[p^{-1}]\rtimes n_0\Z.$
\end{lem}

\begin{proof}
It suffices to note that $xHx^{-1}\cap H$ is either $\Z$ or of the form
$p^n\Z$ for $n\geq 0$, in which case $\sigma(p^nh)=\sigma(h)$ is equivalent to
$p^n\equiv 1\,\operatorname{mod}\,q$.
\end{proof}

Lemma~\ref{B_for_BS}  implies that $\{ngn^{-1}H\mid n\in N,
g\in H\backslash B/H\}$ contains infinitely many disjoint right cosets, and
so we infer the following from \cite[Corollary 1.10]{Bi2} and
Theorem~\ref{commutant_is_Hecke_vNalg}:

\begin{cor} $\mathcal{R}(G, H,\sigma)$ is a factor.
\end{cor}

 Since $B^+=\{x\in B\mid xHx^{-1}\supset H\}= \Z[p^{-1}]\rtimes n_0\N$, the pair
$(B, H)$ is directed. Then $C^*(\H_\sigma(G, H))$ is equal to $p_\sigma C^*(B_\sigma)
p_\sigma$ by Proposition~\ref{C*_completion_B_directed}, and is Morita-Rieffel
equivalent to $C_0(\sigma +H_{\infty}^\perp)\rtimes B/N$ by
Corollary~\ref{B_directed_N_abelian}. Towards describing the last crossed
product we  dwell a little longer on the structure of $N_\sigma=
\Omega_{\mathbf q}$. Note that for $a$ and $b$ as in \eqref{formal_a}
with the sums in $a_{-}$ and $b_{-}$ starting from $-i_0$ and $-j_0$ respectively,
$\frac 1q (a\cdot b)$
is well-defined as an element of $\Q/\Z$, because in the product there are
only finitely many terms not in $\Z$. With $e(x):=\operatorname{exp}(2\pi i
x)$, we claim that
\begin{equation}
\langle a,b\rangle=e(\frac 1q (a\cdot b))\text{ for }a,b\in N_\sigma
\label{duality_qadics}
\end{equation}
is a well-defined duality\footnote{We can also appeal to the second theorem in \cite[\S 12.3.35]{P},
which shows that $\langle a,b\rangle=\operatorname{exp}\Big[2\pi i \sum_{n=-M}^N b_n \Big(\sum_{m=n}^N
\frac {a_m}{{\mathbf q}_n\dots{\mathbf q}_m }\Big)\Big]$ implements a self-duality of $N_\sigma$. The third theorem in
\cite[\S 12.3.35]{P} gives necessary and sufficient conditions on the
${\mathbf q}$-numbers that admit a multiplication, and our choice of
${\mathbf q}$ certainly fulfills those conditions.} pairing.
To prove this claim is essentially an argument similar
to the proof of the second theorem in \cite[\S 12.3.35]{P}, and we
leave the details to the reader.

One can check from \eqref{duality_qadics} that the annihilator of
$\Omega_{\mathbf q}^0$ can be identified as the set of sequences $\{q\sum_0^\infty a_ip^i\}$ with
$0\leq a_i <p$, i.e. as $q\Z_p$.
Since $\sigma(x)=e(\frac 1q x)$, the set defined in \eqref {def_of_H_perp}
of all extensions to $N_\sigma$ is equal to
$$
\sigma + H_\sigma^\perp=\{a\in N_\sigma\mid \langle a,b\rangle=\sigma(b), \forall
b\in \Omega_{\mathbf q}^0\}=1+q\Z_p.
$$

Since $q$ is a unit in $\Z_p$, the element $w_0:=-q^{-1}$ belongs to $\Z_p$. We
let
\begin{equation}
z_0:=1+q w_0 \in \Omega_{\mathbf q}^0\setminus\{0\}.
\label{def_of_z0}
\end{equation}

\begin{lem}
We have $p^{n_0}z_0=z_0$ and $H_\infty=\{jz_0\mid 0\leq j< q\}$,
\cf  Theorem~\ref{theo_p1_B_directed}.
\label{pz0_z0}
\end{lem}

\begin{proof} By the choice of $n_0$, for each
$k\in \N$ there is $s_k\in \Z$ such that $p^{kn_0}=1+qs_k$.
Since $s_1=q^{-1}(p^{n_0}-1)=(1-p^{n_0})w_0$ in $\Z_p$, we have
$$
p^{n_0}z_0=p^{n_0}+qp^{n_0}w_0=1+q(s_1+p^{n_0}w_0)=z_0,
$$
as claimed. For the second claim, we have by definition that $
H_\infty=\bigcap_{k\in \N}p^{kn_0}\Omega_{\mathbf q}^0$,
Clearly $z_0\in H_\infty$, and since $qz_0=q(1+qw_0)=0$ we get $jz_0\in H_\infty$ for
$0\leq j<q$.

To prove the other inclusion, we claim first that $q\Z_p\cap H_\infty=
\{0\}$. Indeed, if $a\in \Z_p$
and $qa\in H_\infty$, then $p^{-n_0}qa\in \Omega_{\mathbf q}^0$, so there are $0\leq b_*<q$ and $b_+\in \Z_p$ such that
\[
qa=p^{n_0}(b_*+qb_+)=(1+qs_1)b_* +qp^{n_0}b_+
=b_*+q(s_1b_*+p^{n_0}b_+).
\]
It follows that $b_*=0$ and $p^{-n_0}a\in \Z_p$. By repeating the argument, we conclude that $a\in \cap_{k\geq 0}p^{kn_0}\Z_p=\{0\}$. To finish off, suppose that
$a=a_*+qa_+$ is in $H_\infty$, where $0\leq a_*<q$ and $a_+\in \Z_p$.
Then $a-a_*z_0\in q\Z_p\cap H_\infty=\{0\}$, and so $a=a_*z_0$, as needed.\end{proof}

\begin{lem}
The annihilator $H_\infty^\perp$ of $H_\infty$ in
$\Omega_{\mathbf q}$ is equal to $\bigcup_{k\geq 0}p^{-kn_0}(q\Z_p)$ and
to the set $Y$ of elements $p^{-kn_0}\alpha_- +\alpha_* +q\alpha_+$ such that
$k>0$, $0\leq \alpha_-<p^{kn_0}$, $\alpha_-+\alpha_*\in q\Z$, and $\alpha_+\in \Z_p$.
\label{Hsigma_perp}
\end{lem}

\begin{proof}
Since by Lemma~\ref{pz0_z0} $H_\infty$ is a cyclic group generated by $z_0$,
$$
H_\infty^\perp=\{z_0\}^\perp=\{a\in N_\sigma\mid az_0\in q\Z_p\}.
$$

Given $a$ in $\Omega_{\mathbf q}$, we can write $a=p^{-kn_0}a_-+a_*+qa_+$ for
$k>0$, $0\leq a_-<p^{kn_0}$, $0\leq a_*<q$, and $a_+\in \Z_p$. By
Lemma~\ref{pz0_z0}, $az_0=a(p^{kn_0}z_0)=(p^{kn_0}a)z_0$, and so
$az_0$ has form $a_-+a_*+qa'_+$ for $a'_+$ in $\Z_p$. Thus $az_0\in q\Z_p$ if
and only if $a_-+a_*\in q\Z$, showing that $\{z_0\}^\perp=Y$.

Since $p^{n_0}z_0=z_0$ and $q\Z_p\subset \{z_0\}^\perp$, we also have
$\bigcup_{k\geq 0}p^{-kn_0}(q\Z_p)\subset \{z_0\}^\perp$. If now $a\in
\{z_0\}^\perp$, we have seen that $a_-+a_*\in q\Z$, and then a calculation shows
that $a\in p^{-kn_0}(q\Z_p)$, finishing the proof.
\end{proof}

\begin{cor}
Let $X_0$ denote the open and  closed subset $(1+q\Z_p)\setminus
p^{n_0}(1+q\Z_p)$ of $\Omega_{\mathbf q}^0$. Then the set $\sigma +H_{\infty}^\perp$
from Corollary~\ref{B_directed_N_abelian} is the disjoint union of $p^{n_0}$-invariant sets
\begin{equation}
\{z_0\}\cup \bigcup_{k\in \Z}p^{kn_0}X_0.
\label{describe_Hinfty}
\end{equation}
\label{H_sigma_infty_perp}
\end{cor}

\begin{proof}
From its definition, the set $\sigma +H_{\infty}^\perp$ is equal to $1+H_\infty^\perp$.
Since
$$
1+p^{-kn_0}(q\Z_p)=p^{-kn_0}(1+qs_k+q\Z_p)=p^{-kn_0}(1+q\Z_p),
$$
Lemma~\ref{Hsigma_perp} implies that $\sigma +H_{\infty}^\perp=
\bigcup_{k\geq 0}p^{-kn_0}(1+q\Z_p)$. The inclusion $p^{ln_0}(1+q\Z_p)\subset 1+q\Z_p$ for all $l\geq 0$
implies that
$$
\bigcup_{k\geq 0}p^{-kn_0}(1+q\Z_p)=\bigcap_{l\geq 0}p^{ln_0}(1+q\Z_p)
\cup \bigcup_{k\in \Z}p^{kn_0}X_0.
$$
To see that this decomposition is exactly \eqref{describe_Hinfty}, we
need to verify that in the right hand side the union is over disjoint sets, and
that
\begin{equation}
\bigcap_{l\geq 0}p^{ln_0}(1+q\Z_p)=\{z_0\}.
\label{z0_intersection}
\end{equation}
First, the inclusions $p^{kn_0}X_0\subset p^{kn_0}(1+q\Z_p)\subset p^{n_0}(1+q\Z_p)$
for
$k> 0$ show that $p^{kn_0}X_0$ and $X_0$ are disjoint, and since $p^{-kn_0}X_0$
is disjoint from $1+q\Z_p$, it is also disjoint from $X_0$. Next, suppose that $z$ is
in the left hand side of \eqref{z0_intersection}. Then
there is $a_l\in \Z_p$ for each $l\geq 0$ such that, with $s_l$ as in the proof of
Lemma~\ref{pz0_z0}, we have
$$
z=p^{ln_0}(1+qa_l)=1+q(s_l+p^{ln_0}a_l).
$$
Since $p^{ln_0}a_l$ is in $p^{ln_0}\Z_p$, it converges to $0\in \Z_p$ as
$l\to \infty$. A computation shows that $s_{l+1}=s_1+s_l+qs_1s_l$ for all
$l\geq 0$, so by passing to a subnet we may assume that $s_l$ converges
to an element $s_*$ that satisfies $s_*=s_1+s_*+qs_1s_*$. This equation has
solution $s_*=-1/q=w_0$ in $\Z_p$, and hence $
z=\lim (1+q(s_l+p^{ln_0}a_l))=1+qw_0=z_0$,
as claimed. To finish, note that if $p^{kn_0}z_0\in X_0$
for some $k\in \Z$, then $z_0\in X_0$, a falsehood.
\end{proof}

The main result concerning the structure of the Hecke algebra is the following.

\begin{thm}
The generalised Hecke $C^*$-algebra $C^*(\H_\sigma(G, H))$ is Morita-Rieffel
equivalent to $C(\T)\oplus (C(X_0)\otimes \mathcal{K}(l^2(\Z))).$
\end{thm}

\begin{proof} By Corollary~\ref{B_directed_N_abelian}, $C^*(\H_\sigma(G, H))$ is
Morita-Rieffel equivalent to $C_0(\sigma +H_{\infty}^\perp)\rtimes B/N$. Applying
Corollary~\ref{H_sigma_infty_perp} gives that $C_0(\sigma +H_{\infty}^\perp)$ is the
direct sum of $C(\{z_0\})$ and $C_0(\bigcup_{k\in \Z}p^{kn_0}X_0)$. But
$\bigcup_{k\in \Z}p^{kn_0}X_0$ is homeomorphic to
$X_0\times \Z$ via a map which is equivariant for the action of $B/N=n_0\Z$,
and the result follows.
\end{proof}

To conclude, we
recollect that for $\H(G, H)$, \cite[Theorem 1.9]{LaLa0}  shows
that the  universal $C^*$-completion $A:=C^*(\H(N\rtimes \Z, \Z))$ is
canonically isomorphic to a semigroup crossed product
$C^*(N/\Z)\rtimes \N$. Then \cite[Theorem 2.1]{BLPR} says that
$A$ has an ideal $C(\mathcal{U}(\Z_p))\otimes
 \mathcal{K}(l^2(\N))$ such that the resulting quotient is $C(\T)$ (here
$\mathcal{U}(\Z_p)$ is the group of units in the ring of $p$-adic integers). The
Schlichting completion of $(G, H)$ is, as noted already, $(\Q_p\rtimes \Z,\Z_p)$.
Then the Morita-Rieffel equivalence implemented by the full projection
$\Chi_{\Z_p}$ carries the ideal $C(\mathcal{U}(\Z_p))\otimes
 \mathcal{K}(l^2(\N))$ of $A$ to the
ideal $C(\mathcal{U}(\Z_p))\otimes \mathcal{K}(l^2(\Z))$
of $C^*(\Q_p\rtimes \Z)$.

\end{ex}

\begin{ex}\label{version_Heisenberg} ($p$-adic version of the Heisenberg group.)
Suppose that $p$ is a prime
and $q,r$ are integers co-prime with $p$ and co-prime with each other. Let
\[
G=\bigl\{\,[u,v,w] \bigm| u,v\in \Z\lbrack p^{-1}\rbrack,
w\in \Z\lbrack p^{-1}\rbrack/\Z\,\bigr\},\]
$H=\bigl\{\, [m,n,0] \bigm| m,n\in \Z\,\bigr\}\cong \Z\times \Z$,
and $\sigma(m, n,0)=\operatorname{exp}(2\pi i (\frac mq+ \frac nr))$.
This is a reduced Hecke triple, and $K=\ker \sigma=q\Z\times r\Z$. It follows as
in examples~\ref{rational_Heisenberg} and \ref{ax_p_adic} that the Schlichting
completion of $(G, H, \sigma)$ consists of
\[
G_\sigma=\bigl\{\,[a,b,w] \bigm| a\in \Omega_{\mathbf q}, b\in
\Omega_{\mathbf r},w\in \Z\lbrack p^{-1}\rbrack/\Z\,\bigr\},
\]
the compact open subgroup $H_\sigma=\bigl\{\,[a,b,0] \bigm| a\in \Omega_{\mathbf q}^0,
b\in \Omega_{\mathbf r}^0,\bigr\}$,
with the natural extension of $\sigma$ to $H_\sigma$.

Since $\sigma$ extends to a
character of $G$, Proposition~\ref{sigma_extends}
implies that $\H_\sigma(G, H)$ is isomorphic to $\H(G, H)$. However, we claim that
$\sigma$ is not continuous for the Hecke topology from $(G, H)$. Indeed,
one can verify that for $x=[p^{-m},p^{-n},0]$ in $G$, where $m,n\geq 0$,
$K\cap xKx^{-1}$ is $p^mq\Z\times p^nr\Z$. The latter set can contain no $H\cap yHy^{-1}$,
which has form $p^k\Z\times p^l\Z$ for $y=[p^{-k}, p^{-l}, w]$ in $G$. Thus the
continuous map $\iota:G_\sigma\to G_0$ is not open.
\end{ex}

\begin{ex}\label{ex:full ax+b} (The full $ax+b$-group of $\Q$.) Let
$N=(\Q,+)$ and consider $Q=(\Q_+^*,\cdot)$ acting by multiplication
$(x,k)\mapsto xk$ for $x\in \Q_+^*,\, k\in \Q$. Then
$G:=N\rtimes Q$ and  $H=\Z$ form the reduced Hecke pair from \cite{BC}.
We can identify $G$ as a matrix group in the form
$$
G=\left\{
\begin{pmatrix}
1&b\\
0&a
\end{pmatrix}
\biggm|\, a\in \Q_+^*, \, b\in \Q\right\},
$$
and then $N$ is the subgroup with $a=1$ in which $H$ is the subgroup with
$b\in \Z$.

Let $n$ be a non-zero positive integer and $\sigma $ the character
of $H$ defined by $\sigma(m)=\exp(2\pi i m/n)$ for $m\in \Z$. Taking
$x_0=(1, n^{-1})$ in $G$ shows that
$K:=\ker \sigma$ contains $x_0Hx_0^{-1}$, and so $\sigma$ is continuous
with respect to the Hecke topology from $(G, H)$.
It is known, see for example \cite{Laca_dil}, that $(G_0:=\mathcal{A}_f\rtimes
Q, H_0:=\mathcal{Z})$ is the Schlichting completion of $(G, H)$, and then
Theorem~\ref{uniqueness_univ_prop} implies that $G_\sigma=G_0$ and $H_\sigma=H_0$.
As a consequence of the description of $B$ below, we see that
 a neighbourhood subbase at $e$ for the Hecke topology from $(B, K)$ consists
of the groups $\{mn\Z\mid m \equiv \pm 1 \bmod n\}$, and so does not
contain the open set $n^2\Z$ in $G_0$.
It was computed in \cite[Lemma 3.2.3]{Cu1} that $B$ consists of pairs
$(r, a/b)$ with $r\in \Q, a,b\in \N, b>0$, and such that
$\text{gcd}(a, b)=1$, and $a-b\in n\Z$. We claim (without proof) that the following
description of $B$ is valid:

\begin{lem} Let $T_q$ be the subsemigroup $\{m\in \N^*\mid m-1\in n\N\text{ or }
m+1\in n\N\}$ of $\N^*$. Then $B$ is the subgroup $\Q\rtimes T_qT_q^{-1}$ of $G$,
and $(B, H)$ is directed.
\end{lem}

Here $G_\sigma$ is the same for all $\sigma$ and we want to study the relation
between the different ideals $\overline{C^*(G_\sigma)p_\sigma C^*(G_\sigma)}$.
By Theorem~\ref{Hecke_alg_as_cp} this is the same as studying the ideals
$I_{\sigma}$ in $C^*(N_\sigma)$, or by Theorem~\ref{theo_N_abelian} the sets
$\Omega_\sigma$ defined in \eqref{def_omega}. If $\sigma$ corresponds
to $n\in \N^*$ then, since $N_\sigma=\mathcal{A}_f$ and
$\mathcal{A}_f$ is self-dual with duality carrying $\mathcal{Z}$ into
$\mathcal{Z}^\perp$, we have that
$
\sigma +H_\sigma^\perp=\{a\in \mathcal{A}_f\mid a-\sigma\in \mathcal{Z}^\perp\}=1/n +\mathcal{Z}$.
Thus the sets are given by
$$
\Omega_n=\bigcup_{t\in Q}t(\frac 1n +\mathcal{Z}).
\label{Omega_q_ax+b_group}
$$
We can then describe exactly the ideals $C_0(\Omega_n) \rtimes
Q$ inside $C_0(\mathcal{A}_f)\rtimes Q$ corresponding to different choices of
$\sigma$, and link to the results of \cite{LaR2} and \cite{BLPR}.

\begin{lem}
\textnormal{(a)} If $q$ is a prime and $m>0$, then $\Omega_{q^m}=\{x\in \mathcal{A}_f\mid x_q\neq 0\}$.
In particular, $\Omega_{q^m}$ is independent of $m$.

\smallskip
\textnormal{(b)} If $n=q_1^{i_1},\cdots ,q_l^{i_l}$ with the $q_j$'s
different primes and each $i_j>0$, then $$ \Omega_{n}=\{x\in
\mathcal{A}_f\mid x_{q_j}\neq 0,\, j=1,\cdots ,l\}.
$$
\label{lemma_Omega}
\end{lem}

\begin{proof} We only prove part (b) for $l=2$, as the rest follows by similar
arguments. Thus we assume that $n=q^j r^i$ with $q$ and $r$
distinct primes.

The forward inclusion is obvious. To prove the other,
we take $x\in \mathcal{A}_f$ with $x_q\neq 0$
and $x_r\neq 0$.
We can write $x_q=q^l(y+q^j y_q)$ and $x_r=r^m(w+r^i w_r)$ with
$0<y<q^j$, $q$ not dividing $y$, $y_q\in \Z_q$, $0<w<r^i$, $r$ not
dividing $w$,
 and $w_r\in \Z_r$. Pick $s,t,c,d\in \Z$ such that
$sy- tq^j =1$ and $cw-d r^i=1$. By the Chinese Remainder Theorem there is
an integer
$k$ such that
$$
k\equiv sr^m\,\bmod{q^j}\text{ and }k\equiv cq^l\,\bmod{r^i}.
$$
Then $ky\equiv r^m\,\bmod{q^j}$ and $kw\equiv q^l\,\bmod{r^i}$, and so
there are $y'_q\in \Z_q$ and $w'_r\in \Z_r$ such that
$$
x_q=q^lr^mk^{-1}(1+q^j r^i y'_q)\text{ and }x_r=q^lr^mk^{-1}(1+q^j r^i
w'_r).
$$
Since $q$ and $r$ are units in $\Z_p$ for every prime $p$ different from
$q$ and
$r$, there is $z_p\in \Z_p$  such that $x_p=q^lr^mk^{-1}(1+ q^j r^i z_p)$.
Hence
with $z:=(z_p)\in \mathcal{Z}$, and replacing (if necessary) $k$ with $-k$,
we have $x\in \Omega_{n}$, as claimed.
\end{proof}
\end{ex}

We can now resume our description of $\H_\sigma(G, H)$. When $p$ is a prime,
\cite[Proposition 2.5]{LaR2} says that $
J_p:=C_0(\mathcal{A}_f\setminus \{x\in \mathcal{A}_f\mid x_p=0\})\rtimes Q
$ is one of the primitive ideals of the dilation
$C_0(\mathcal{A}_f)\rtimes Q$ cf. \cite{Laca_dil} of the Hecke $C^*$-algebra of Bost and Connes,
see also \cite[\S4]{BLPR}. Suppose that $n$ is a non-zero positive integer and
let $S$ be the set of primes in the decomposition of $n$. Using
Lemma~\ref{lemma_Omega} shows that $\H_\sigma(G, H)$ is Morita-Rieffel
equivalent to
$$
C_0(\Omega_n)\rtimes Q=C_0(\bigcap_{p\in S} \Omega_p)\rtimes Q=\bigcap_{p\in S}J_p.
$$

\begin{ex}\label{ex:lamp} (The lamplighter group, see e.g. \cite{dlH}.)
{\rm
Suppose that $F$ is a finite abelian group with identity $e$. Let
$$
N_{-}=\bigoplus_{-\infty}^0
F,\quad H=\bigoplus_{1}^\infty F,
$$
and set $N=N_{-}\oplus H$. The forward shift $\alpha$ on $N$ acts as
$\alpha((x_k)_{k\in \Z})=(y_k)_{k\in
\Z}$, with $y_k=x_{k-1}$ for all $k\in \Z$.

It is proved in \cite[Lemma 3.1.1]{Cu1} that $(G:=N\rtimes_\alpha \Z, H)$
is a Hecke pair,  and we note that
$(G, H)$ is reduced. Let $H_0$ be the profinite, hence compact, abelian group
$\varprojlim_{n\geq 1}F^n$, identified as $\prod_{n=1}^\infty
F$, and set $N_0:=N_{-}\oplus H_0$. We regard an element of $N_0$ as a sequence
$(x_k)_{k\in \Z}$ in $\prod_{k\in Z}F$ such that for some integer $n_0$ we have
$x_k=0$ for $k<n_0$.  Let $\beta$ be the natural continuous extension of
$\alpha$ to $N_0$.
The inclusion of $H$ in $H_0$ is equivariant for the actions of $\Z$, and so gives rise to a homomorphism
$\phi:N\rtimes_\alpha \Z\to N_0\rtimes_\beta \Z$.
By construction, $\phi$ has dense range and $\phi^{-1}(H_0)=H$. Hence \cite[Theorem 3.8]{KLQ} implies that $(G_0:=N_0\rtimes_\beta \Z, H_0)$
is the Schlichting completion of $(G, H)$.

The subset $T=\{(y,k)\in G\mid (y,k)H(y,k)^{-1}\supseteq H\}$ is a
subsemigroup of $G$, and
$T=\{(y,k)\in G\mid k\geq 0\}$. Thus $G=T^{-1}T$, so $(G, H)$ is directed in
the sense of \cite[\S 6]{KLQ}. By \cite[Theorems 6.4 and 6.5]{KLQ},  or by
\cite[Theorem 1.9]{LaLa0}, the universal $C^*$-completion
$C^*(G, H)$ of $\H(G, H)$ is isomorphic to the
corner in $C^*(G_0)$ determined by the full projection $\Chi_{H_0}$. Moreover,
\cite[Corollary 8.3]{KLQ} (or \cite[Theorem 2.5]{LarR}) imply that $C^*(G, H)$ is
Morita-Rieffel equivalent to $C_0(\widehat{N_0})\rtimes_{\widehat{\beta}} \Z$, where
$
\widehat{N_0}=\bigcup_{n\in \Z}\widehat{\beta_n}(H_0^\perp).$

We now assume that $\sigma$ is a character of $H$.
Thus $\sigma=(\sigma_n)_{n\geq 1}$,
where $\sigma_n$ is a character of $F$ for every $n\geq 1$.
Note that $\sigma(H)$ is included in the finite set
$\{\langle \pi, f\rangle_F\mid \pi\in \widehat{F}, f\in F\}$, where
$(f,f')\mapsto \langle f, f'\rangle_F$ is a fixed self-duality of $F$. However,
$(G, K)$ will not be directed for arbitrary $\sigma$, and to proceed we
specialise further.

When $\sigma$ is not periodic, $B=N$ by \cite{Cu1}. Hence
\cite[Proposition 1.5.2]{Cu1} or
Lemma~\ref{H_sigma_spanned_by_epsilonx} imply that  $\H_\sigma(G, H)\cong
\mathbb{C}(N/H)$. Therefore $C^*(\H_\sigma(G, H))$  is the group algebra $C^*(N_{-})$.

Next we restrict the attention
to $1$-periodic
characters, so we assume that $\sigma_n=\sigma_1$, $\forall n\geq 2$. We have that
$K=\ker \sigma=\{f=(f_k)_k\in H\mid \sigma_1(\sum_k f_k)=1\}$, and $(G, K)$ is reduced.

Let $M$ be the profinite group $\varprojlim_{n\geq 1}(F^n\oplus \sigma(H))$ with
bonding homomorphisms
$$
(f_1,\dots, f_n, f_{n+1},\sigma(h))\mapsto
(f_1, \dots, f_{n}, \sigma(h))
$$
from $F^{n+1}\oplus \sigma(H)$ onto $F^n\oplus \sigma(H)$ and canonical homomorphisms
$\pi^n$  from $M$ onto $F^n\oplus \sigma(H)$. By viewing  $M$ as the subset of
$\prod_{n\geq 1}F^n\oplus \sigma(H)$
of sequences compatible under all bonding maps, we see that for $n\geq 1$ and
$h=(h_k)_{k\geq 0}\in H$, the formula
$
\pi_n(h):=(h_1, \dots ,h_n, \sigma(h))$
defines an element of $M$.  This gives a  homomorphism $\pi:H\to M$ such that $
\pi(h):=(\pi_n(h))_{n\geq 1}\text{ for }h\in H$. Since the
cylinder sets $\{x\mid \pi^n(x)=\pi_n(h)\}$ form a basis for the topology on $M$,
$\pi$ has dense range. Let $L$ be the locally compact group $
L:=(N_{-}\oplus \prod_{n\geq 1}F\oplus \sigma(H))\rtimes_\beta \Z,
$
with the compact, open subgroup $M$, and define a
homomorphism $\phi:G\to L$  and a continuous character $\rho:M\to \T$ by
$
\phi(y,h,k)=(y, \pi(h),k)$,
and $\rho(x, s)=s$ for $(x,s)\in H_0\oplus \sigma(H)$. The map $\phi$ has dense range because $\pi$ does,
and clearly $\phi^{-1}(M)=H$. For $h=(h_1,\dots,h_n, e,\dots)$ in $H$ we have
$$
\rho\circ \phi(h)=\rho\circ\pi(h)=
\rho(h,\sigma(h))=\sigma(h).
$$
Hence Theorem~\ref{uniqueness_univ_prop} yields the following result, which
includes \cite[Example 4]{Cu0} (or \cite[\S 3.1]{Cu1}).

\begin{cor}
The triple $(L, M, \rho)$ is the Schlichting completion $(G_\sigma, H_\sigma, \sigma)$  of $(G, H,\sigma)$. Moreover, the closure of $N$ in the
Hecke topology from $(G, K)$ is
$$
{N}_\sigma=N_{-}\oplus {H}_\sigma.
$$
\label{one_periodic}
\end{cor}

So $G_\sigma$ and $G_0$ are different in this example.
Since $B=G$ by \cite[\S 3.1]{Cu1} and $H\trianglelefteq N\trianglelefteq G$,
Corollary~\ref{C*_completion_normal_subgroups} implies that
$p_\sigma C^*(G_\sigma) p_\sigma$ is the enveloping $C^*$-algebra of
$\H_\sigma(G, H)$. By Theorem~\ref{theo_N_abelian}, this $C^*$-completion is
Morita-Rieffel equivalent to $C_0(\Omega_\sigma)\rtimes_{\beta}\Z$, where
$\Omega_\sigma$ is defined in \eqref{def_omega}. Using results from \cite[\S 3.1]{Cu1}, we will identify this set. The self-duality $(f,f')\mapsto \langle
f, f'\rangle_F$ of $F$ implements a self-duality
$$
\langle (y_k)_{k\in \Z}, (y_k')_{k\in \Z}\rangle
=\prod_{n\geq 1}\langle y_{1-n},y_n'\rangle_F \prod_{n\geq 1}\langle
y_{1-n}', y_n\rangle_F
$$
of $N_0=N_{-}\oplus H_0$. Note that under this identification
${H_0}^\perp$ is carried into $H_0$. The group $\sigma(H)$ is
also self-dual with duality expressed in terms of a fixed generator $\omega$
as $\langle \omega^j, \omega^k\rangle=\omega^{jk}$ for $j, k\in \Z$. Hence
${N}_\sigma=N_0\oplus \sigma(H)$ is
self-dual. Under the described duality pairings, an element $(x, \omega^j)$ of
${N}_\sigma$ is in $\sigma +H_\sigma^\perp$ precisely when $j=1$ and
$x\in {H_0}^\perp$. Thus $\sigma +H_\sigma^\perp={H}_0\oplus\{\omega\}$.

The action of $\Z$ on $N_\sigma$ is $\beta$ on $N_0$ and is trivial
on $\sigma(H)$. The formulas in
\cite[Lemma 3.1.6]{Cu1} show that the same is true for the dual action $\hat{\beta}$,
and then
$$
\Omega_\sigma=\bigcup_{n\in \Z}\hat{\beta}_n (\sigma +H_\sigma^\perp)
=(\bigcup_{n\in \Z}\hat{\beta}_n (H_0))\oplus \{\omega\}=N_0\oplus\{\omega\}.
$$
Finally, we get  an isomorphism
$
C_0(\Omega_\sigma)\rtimes_{\hat{\beta}} \Z\cong (C_0(N_0)\rtimes_{\hat{\beta}} \Z)
\oplus C(\T)$;
in other words, the ideal $\overline{C^*(G_\sigma)p_\sigma C^*(G_\sigma)}$
is isomorphic to $C^*(G_0)\oplus C(\T)$.
}\end{ex}

\end{document}